%% file: main.tex
\documentclass{article}

\usepackage{PRIMEarxiv}

\usepackage[utf8]{inputenc} % allow utf-8 input
\usepackage[T1]{fontenc}    % use 8-bit T1 fonts
\usepackage{hyperref}       % hyperlinks
\usepackage{url}            % simple URL typesetting
\usepackage{booktabs}       % professional-quality tables
\usepackage{amsfonts}       % blackboard math symbols
\usepackage{nicefrac}       % compact symbols for 1/2, etc.
\usepackage{microtype}      % microtypography
\usepackage{lipsum}
\usepackage{fancyhdr}       % header
\usepackage{graphicx}       % graphics
\graphicspath{{figures/}}     % organize your images and other figures under media/ folder
\usepackage{amsmath}
\usepackage{amsthm}
\usepackage{amssymb}
\usepackage{multirow}
\usepackage{amsfonts}
\usepackage{bm}
\usepackage{siunitx}
\usepackage{tabularx} 
\usepackage{multicol} 
\usepackage[table]{xcolor}
\usepackage[ruled,vlined]{algorithm2e}
\usepackage{comment}

\usepackage{etoolbox}
\usepackage{nomencl}
\makenomenclature

\renewcommand{\nomgroup}[1]{%
  \ifstrequal{#1}{L}{\item[\textbf{Latin Symbols}]}{%
  \ifstrequal{#1}{G}{\item[\textbf{Greek Letters}]}{%
  \ifstrequal{#1}{A}{\item[\textbf{Acronyms}]}{}}}}

\title{A Comparison of Parametric Dynamic Mode Decomposition Algorithms for Thermal-Hydraulics Applications}

\author{
  Stefano Riva, Andrea Missaglia, Carolina Introini \\
  Department of Energy\\
  Politecnico di Milano\\
  Milan, Italy 20156 \\
  \texttt{stefano.riva@polimi.it; andrea.missaglia@polimi.it; carolina.introini@polimi.it} \\
  \AND
   In Cheol Bang \\
   Nuclear Thermal-hydraulics and Safety Lab\\
   Ulsan National Institute of Science and Technology\\
   Ulsan,  Republic of Korea \\
    \texttt{icbang@unist.ac.kr} \\
    \AND
    Antonio Cammi \\
    Emirates Nuclear Technology Center, Department of Mechanical and Nuclear Engineering \\
    Khalifa University\\
    Abu Dhabi, United Arab Emirates \\
     (and) Department of Energy\\
     Politecnico di Milano\\
     Milan, Italy 20156 \\
     \texttt{antonio.cammi@ku.ac.ae} \\
}

\DeclareMathOperator*{\argmin}{arg\,min}

\renewcommand{\vec}[1]{\mathbf{#1}}

\newcommand{\norma}[1]{\left|\left|{#1}\right|\right|}

\newcommand{\dpart}[2]{\frac{\partial {#1}}{\partial {#2}}}
\newcommand{\der}[2]{\frac{d {#1}}{d {#2}}}

\begin{document}
\maketitle

\begin{abstract}
In recent years, algorithms aiming at learning models from available data have become quite popular due to two factors: 1) the significant developments in Artificial Intelligence techniques and 2) the availability of large amounts of data. Nevertheless, this topic has already been addressed by methodologies belonging to the Reduced Order Modelling framework, of which perhaps the most famous equation-free technique is Dynamic Mode Decomposition. This algorithm aims to learn the best linear model that represents the physical phenomena described by a time series dataset: its output is a \textit{best-fit} linear operator of the underlying dynamical system that can be used, in principle, to advance the system itself in time even beyond the training time interval. However, in its standard formulation, this technique cannot deal with parametric time series, meaning that a different linear model has to be derived for each configuration. Research on this is ongoing, and some versions of a \textit{parametric Dynamic Mode Decomposition} already exist. This work contributes to this research field by comparing the different algorithms presently deployed and assessing their advantages and shortcomings compared to each other. To this aim, three different thermal-hydraulics problems are considered: two benchmark 'flow over cylinder' test cases at diverse Reynolds numbers, and the DYNASTY experimental facility operating at Politecnico di Milano, which studies the natural circulation established by internally heated fluids for Generation IV nuclear applications, simulated using the RELAP5 nodal solver. As a key result, this paper highlights the main advantages and disadvantages of the available parametric Dynamic Mode Decomposition methods, concluding that the choice of the algorithm version strongly depends on the problem under consideration and on the user priority (i.e., accuracy or computational speed of the online phase); additionally, this paper shows that an interpretable linear model can be learnt from parametric datasets governed by non-linear models, that can be used for parameter interpolation and, most importantly, for state prediction in time.
\end{abstract}

% keywords can be removed

\section{Introduction}
In the broad world of computational sciences for engineering applications, there is a need for accurate, reliable and efficient mathematical models for multi-query or real-time scenarios \cite{quarteroni_reduced_2015, rozza_model_2020}, both for design and monitoring purposes \cite{mohanty_development_2021}. Nowadays, mathematical modelling is (mainly) carried out through Partial Differential Equations (PDEs), and their discretized version (retrieved adopting the most suitable technique for the application of interest) represents the Full Order Model (FOM), which typically describes, typically in detail, the underlying physics of the investigated system. The numerical solution of this FOM is, however, often computationally expensive, which makes it not suited for all these scenarios requiring multiple in-series solutions (design, parameter optimization, sensitivity analyses) or real-time evaluations of the state of the system (online monitoring, control) \cite{gong_data-enabled_2022}: to address this issue, the Reduced Order Modelling (ROM) framework \cite{rozza_model_2020, lassila_model_2014, brunton_data-driven_2022} has been introduced. These methodologies allow retrieving a surrogate (i.e., low dimensional) mathematical model starting from a high-fidelity and highly complex one, to be used for multi-query applications, where multiple evaluations of the state of the dynamical system under different conditions are required, and online monitoring and control, as they greatly reduce the computational time and resources needed to obtain a solution through simulation (to the point that such models, once deployed, can even run on personal laptops). Indeed, in a nutshell, all ROM techniques aim at discovering a reduced representation of the FOM starting from some training solutions, called \textit{snapshots}, which accurately describe the physical problem. This reduction becomes functional if the computational time to obtain an approximation of the FOM solution for an unseen parameter (including time) is much shorter compared to that of the FOM, keeping, at the same time, the required accuracy. 

Various ROM techniques exist: among them, the Dynamic Mode Decomposition (DMD) \cite{schmid_dynamic_2010, j_nathan_kutz_chapter_nodate-2} has gained popularity in the engineering community due to 1) its simple implementation and 2) its good performance for different applications, which range from fluid dynamics \cite{hess_data-driven_2023} to plasma physics \cite{faraji_data-driven_2024-1, faraji_data-driven_2024-2} to reactor dynamics \cite{GONG2020107826, INTROINI2024113477, bopdmd_Haidy}. However, the biggest selling point of DMD is that it identifies the dominant structures in the original problem, with which it generates a surrogate linear model based on the starting dataset only, without requiring any knowledge of the governing equations: by retrieving a surrogate linear operator that describes the dataset, it can be used also to advance the dynamical system in time beyond the span of the starting dataset. In fact, DMD defines dimensionality reduction in a pure data-driven (non-intrusive) way, by learning the best linear model underlying the physical phenomena described by the available spatio-temporal data. As said, DMD methods yield to interpretable dynamics which many machine learning models, typically black-boxes, lack \cite{Guo2025_DMDcomparison, DING2024123138, HE2024107258}; furthermore, the amount of data required for training is much less with respect to other surrogate modelling methods for time-series data, like DeepONet or foundation models \cite{LU2025114184, WANG2025111427}. In the nuclear engineering community, it is important to understand how a certain prediction from the mathematical model has been obtained, which is why the choice of DMD over other machine learning methods. However, in its standard formulation, DMD cannot deal with parametric time series, meaning that a different linear model has
to be derived for each parameter realisation. 

Indeed, many systems depend not only on time but also on some parameters, such as the inlet fluid velocity in fluid dynamics, the diffusion coefficient in a nuclear reactor or the ambient temperature in a natural circulation loop. This parametric dependence is relevant in the design phase, where the effect of different parameter values on the operating conditions of the system should be studied, but also during online monitoring and control, as changes from the nominal values may bring the system towards accidental conditions. Thus, both time dynamics and parametric dependence should be learnt to have an 'as-general-as-possible' surrogate model that can describe the system even for unseen scenarios \cite{GONG2020107826}. Parametric ROM techniques are currently under study and development, including parametrised versions of DMD.

To the best of the authors’ knowledge, the most recent and more in-depth research on this topic can be found in \cite{faraji_data-driven_2024-2, andreuzzi_dynamic_2023, huhn_parametric_2023}. Each one of these research works presents different versions of a Parametric Dynamic Mode Decomposition (pDMD) algorithm, applied to several parametric dynamical systems. Andreuzzi et al. \cite{andreuzzi_dynamic_2023} first performs dimensionality reduction of the parametric dataset; then, they feed the reduced dataset for each parameter realization to the DMD algorithm, forecasting the reduced snapshots in the future and then reconstructing the parametric dynamics during the online phase using an approach similar to POD with Interpolation \cite{demo_complete_2019} is used to estimate the parametric dynamics during the online phase. Huhn et al. \cite{huhn_parametric_2023}, instead, collect the eigenvalues and the eigenmodes computed by DMD for each parametric dataset, and interpolate them for the unseen parameters; Faraji et al. \cite{faraji_data-driven_2024-2} extends the previous work to the optimised DMD algorithm \cite{optdmd2018, bopdmd_2022}. 

In \cite{NUTHOS24_pDMD}, the authors have proposed another strategy which builds from the approach described in \cite{andreuzzi_dynamic_2023}, focusing on the estimation of the generic state matrix of the linearised parametric dynamical system. Compared to \cite{andreuzzi_dynamic_2023}, the proposed method avoids the time-advancement through DMD for each parameter realisation in the training set and the online generation of the interpolation maps by restricting the latter in the training phase to keep the online step as fast as possible. Compared to \cite{huhn_parametric_2023}, it avoids interpolating the full-order spatial eigenmodes, which requires high-dimensional interpolants which are highly inefficient from the numerical point of view: instead, the state matrices of each DMD model, computed for each parametric time-series realization of the training dataset, are unfolded and collected as a new series of snapshots, upon which a Singular Value Decomposition (SVD) is performed to retrieve the hidden patterns in the parametric dynamics themselves. All strategies reported here are purely data-driven, in the sense that the temporal dynamics are learnt by the DMD itself and the parametric dependence of the generalised reduced system is obtained through regression techniques. Variants of the DMD in which physical constraints can be enforced are available in the literature \cite{Baddoo2023}, but they are out of the scope of this paper.

This work investigates the three DMD approaches cited in the previous paragraph by applying them to the classical thermal-hydraulics benchmark case of the (laminar) 'flow over cylinder'. In particular, two different datasets have been considered: for the first one, the case study was solved using the \textit{dolfinx} interface for FEniCS, a finite element solver \cite{alnaes_unified_2014, baratta_dolfinx_2023, scroggs_basix_2022, scroggs_construction_2022}; for the second one, the dataset was taken from the CFDbench dataset \cite{luo_cfdbench_2024}. The key difference between the two datasets lies in the parameter range for the Reynolds number: the FEniCS case spans a wider range, whereas the CFDbench one focuses on a smaller range; additionally, slight differences can be found in the geometry and in the inlet boundary conditions for velocity. The chief goal of this comparison is to highlight the advantages and shortcomings of each version of parametric DMD, evaluating their performance in predicting the dynamics of parametric flow data. Following this comparison, the pDMD algorithms will be applied to the study of the DYNASTY experimental facility \cite{cammi_dynasty_2016, pini_experimental_2016, battistini_development_2021, benzoni_preliminary_2023} deployed at Politecnico di Milano, which studies the dynamics of internally heated fluids for Generation-IV reactors: for this case, the FOM is represented by a RELAP5 \cite{fletcher1992relap5} model. Additionally, for this last case, the pDMD prediction has been compared with the available experimental data from the facility. The overall goal of this comparison is to determine which one of the available algorithms is the most suited for online monitoring (and possibly control) for engineering systems. This work acts as a first comparison of different parametric DMD models with applications to thermal-hydraulics systems, focusing on the advantages and disadvantages of the different approaches more than on the complexity of the case studies. Indeed, this study represents a fundamental step before applying these algorithms to more complex situations such as transient studies in nuclear reactors, where more complex dynamics are involved. It is worth mentioning that, in principle, any non-linear dynamical system can be described with a linear dynamics, according to Koopman theory \cite{brunton_data-driven_2022}, provided that a suitable reduced coordinate system has been selected. For more complex, strongly non-linear test cases, more advanced reduction techniques may be used instead of the SVD; however, it is also worth mentioning that, by design, nuclear reactor transients are, in general, characterised by smooth-enough dynamics, and not by strongly non-linear and chaotic phenomena.

The paper is structured as follows: Section \ref{sec: basic-dmd} provides a brief overview of the Dynamic Mode Decomposition, including its basic and optimized versions; Section \ref{sec: pdmd} focuses on the extension to parametric dataset, briefly presenting each algorithm version that will be used within the paper; Section \ref{sec: num-res} provides the key numerical results of the algorithms for the three considered cases: two laminar flow over a cylinder with different parameter ranges and the DYNASTY experimental facility; finally, Section \ref{sec:conclusion} concludes the paper and discusses potential future works and extensions.

\section{Basics of Dynamic Mode Decomposition}\label{sec: basic-dmd}

In this Section, the basics of the Dynamic Mode Decomposition method are briefly recalled. This algorithm was developed by Schmid in \cite{schmid_dynamic_2010} to identify spatial and temporal coherent structures from high-dimensional data for fluid-dynamics applications. This method falls in the broad category of Reduced Basis methods \cite{quarteroni_reduced_2015}, whose rationale is briefly described in Appendix \ref{app: rb}. In particular, DMD is based on the Singular Value Decomposition technique \cite{brunton_data-driven_2022, rozza_model_2020}, and it provides a modal decomposition where each mode consists of spatially correlated structures with the same linear behavior in time (e.g., oscillations at a given frequency with growth or decay); thus, DMD finds the linear coordinate system that best fits the starting dataset. By ranking the importance of the modes through the singular values of the data matrix, which offer a measure of the 'information content' of the associated modes, DMD can build a low-dimensional surrogate linear model that can simulate the temporal evolution of these spatial modes. In the following, the original algorithm by Schmid \cite{schmid_dynamic_2010} is described, then the optimised version, able to better handle non-linearities \cite{optdmd2018} is presented.

\subsection{Original algorithm}
Let $\vec{x}(t)\in\mathbb{R}^{\mathcal{N}_h}$ being some high-dimensional measurements of the state of the system over time $t$, collected at different time instants $\{t_k\}_{k=1}^{N_t}$ and arranged in a snapshot matrix of size $\mathcal{N}_h \times \mathcal{N}_t$ in the following way:
\begin{equation}
    \mathbb{X} = [\mathbf{x}(t_1) \mid \mathbf{x}(t_2) \mid \dots \mid \mathbf{x}(t_{N_t})] =  [\mathbf{x}_1 \mid \mathbf{x}_2 \mid \dots \mid \mathbf{x}_{N_t}] \in \mathbb{R}^{\mathcal{N}_h \times N_t},
\end{equation}
where $\mathcal{N}_h$ is the number of spatial degrees of freedom of the FOM and $\mathcal{N}_t$ is the number of saved time instants. The DMD seeks the best-fit linear operator $\mathbb{A}$ that best represents the dataset $\mathbb{X}$ and that allows the underlying system to advance in time, such that
\begin{equation}
    \mathbf{x}_{k+1} \approx \mathbb{A}\mathbf{x}_k 
\end{equation}
Thus, the operator $\mathbb{A}$ describes the dynamics of the associated \textbf{linear} system that best advances the snapshots forward in time; in control theory, this operator is commonly called state matrix \cite{Nise2000}, and for the rest of the work, the two terms will be used interchangeably. Mathematically, this operator is the solution to the following minimisation problem \cite{faraji_data-driven_2024-1, Tu2014}:
\begin{equation}
    \mathbb{A} = \argmin\limits_{\mathbb{A^\star}\in\mathbb{R}^{\mathcal{N}_h\times \mathcal{N}_h}} \| \mathbb{X}^+ - \mathbb{A}^\star\mathbb{X}^-\|_F = \mathbb{X}^+(\mathbb{X}^-)^\dagger
\end{equation}
given the snapshots matrices $\mathbb{X}^- = [\mathbf{x}(t_1) \mid \dots \mid \mathbf{x}(t_{N_t - 1})] \in \mathbb{R}^{\mathcal{N}_h \times (N_t - 1)}$ and $\mathbb{X}^+ = [\mathbf{x}(t_2) \mid \dots \mid \mathbf{x}(t_{N_t})] \in \mathbb{R}^{\mathcal{N}_h \times (N_t - 1)}$, $\|\cdot \|_F$ indicating the Frobenius norm and the superscript $^\dagger$ indicating the Moore-Penrose pseudo-inverse \cite{brunton_data-driven_2022}. This procedure is known as exact DMD, since the high-dimensional matrices are used. However, it becomes unfeasible when the dimensionality of the data becomes large, i.e. $\mathcal{N}_h >> 1$. Schmid \cite{schmid_dynamic_2010} proposed a reduced version of DMD, where the data matrices are reduced using the Singular Value Decomposition or the equivalent Proper Orthogonal Decomposition (POD) \cite{rozza_model_2020, brunton_data-driven_2022} to retrieve the reduced dimension that guarantees suitable accuracy. The operator $\mathbb{A}$ is thus projected onto the POD modes, and the best-fit operator is now learnt in the reduced space spanned by the POD modes themselves. Briefly, the SVD decomposes a generic matrix $\mathbb{X}$ into three contributions\footnote{The matrix $\mathbb{X}$ is assumed to be real, which is a suitable assumption for nuclear reactors data and more in general for engineering applications. This theory can be extended to complex matrices naturally.} in the following way:
\begin{equation}
    \mathbb{X} = \mathbb{U}\Sigma\mathbb{V}^T
\end{equation}
where $\mathbb{U}$ encodes the dominant spatial structures, named \textbf{modes}, $\Sigma$ is a diagonal matrix containing the singular values associated with the modes in $\mathbb{U}$ and $\mathbb{V}$ embeds the temporal dynamics. In practice, this decomposition implies that any matrix can be approximated with a low-rank representation: in fact, SVD sorts the modes hierarchically according to their content of energy, that is, how much information of the original dataset they retain \cite{rozza_model_2020, brunton_data-driven_2022}; the singular values associated with each mode indicate this information quantity. It is, therefore, possible to truncate the SVD matrices, and thus retrieve a reduced representation of the original dataset, to a rank $r$\footnote{The rank $r$ is a hyperparameter that must be tuned according to the desired level of accuracy, but also according to the application. A possible way to calibrate this consists in looking at the decay of the singular values, computed by the SVD \cite{quarteroni_reduced_2015, brunton_data-driven_2022}.} by retaining only the first $r$ terms as follows:
\begin{equation*}
    \mathbb{X}\approx \widetilde{\mathbb{U}}\widetilde{\Sigma}\widetilde{\mathbb{V}}^T
\end{equation*}
where $\widetilde{\mathbb{U}}\in\mathbb{R}^{\mathcal{N}_h\times r}$, $\widetilde{\Sigma}\in\mathbb{R}^{r_h\times r}$ and $\widetilde{\mathbb{V}}\in\mathbb{R}^{N_t\times r}$. Adopting this approximation, it is possible to obtain a more computationally efficient DMD and to compute the low rank operator $\widetilde{\mathbb{A}}\in\mathbb{R}^{r\times r}$ as the projection of the operator $\mathbb{A}$ onto the first $r$ SVD modes of matrix $\mathbb{X}^-$
\begin{equation}
    \widetilde{\mathbb{A}} = \widetilde{\mathbb{U}}^T\widetilde{\mathbb{A}}\widetilde{\mathbb{U}}= \widetilde{\mathbb{U}}^T\mathbb{X}^+\widetilde{\mathbb{V}}\widetilde{\Sigma}^{-1}
\end{equation}
Tu et al. \cite{Tu2014} proved the strong connection between the reduced eigenvalues of the low-dimensional operator $\widetilde{\mathbb{A}}$ and the true ones of $\mathbb{A}$, so that the former can be used in place of the latter. The DMD algorithm leverages this dimensionality reduction to evaluate the dominant eigenvalues and eigenvectors of the best-fit linear operator, thus avoiding the computation of the pseudo-inverse $(\mathbb{X}^-)^\dagger$. The overall procedure to discover this operator can be found in Algorithm 1.

\begin{algorithm}[htbp]
\caption{Dynamic Mode Decomposition: operator learning}
\KwIn{Snapshot Matrix $\mathbb{X} = [\mathbf{x}(t_1) \mid \mathbf{x}(t_2) \mid \dots \mid \mathbf{x}(t_{N_t})] \in \mathbb{R}^{\mathcal{N}_h \times N_t}$}
\KwOut{Best fit operators $\mathbb{A}$ and $\widetilde{\mathbb{A}}$}
\vspace{0.25cm}
\textbf{Generate time advancement matrices:} \\
$\mathbb{X}^- = [\mathbf{x}(t_1) \mid \dots \mid \mathbf{x}(t_{N_t - 1})] \in \mathbb{R}^{\mathcal{N}_h \times (N_t - 1)}$ \;
$\mathbb{X}^+ = [\mathbf{x}(t_2) \mid \dots \mid \mathbf{x}(t_{N_t})] \in \mathbb{R}^{\mathcal{N}_h \times (N_t - 1)}$ \;

\textbf{Compute the Singular Value Decomposition (SVD):} $\mathbb{X}^- = \widetilde{\mathbb{U}} \widetilde{\Sigma} \widetilde{\mathbb{V}}^T$ \;

\textbf{Compute best fit operator:} $\mathbb{A} = \mathbb{X}^+ \widetilde{\mathbb{V}} \widetilde{\Sigma}^{-1} \widetilde{\mathbb{U}}^* \in \mathbb{R}^{\mathcal{N}_h \times \mathcal{N}_h}$ \;

\textbf{Project} $\mathbb{A}$ \textbf{onto the modes with rank} $r$: \\$\widetilde{\mathbb{A}} = \widetilde{\mathbb{U}}^* \mathbb{A} \widetilde{\mathbb{U}} \in \mathbb{R}^{r \times r}$
\end{algorithm}

Typically, the advancement in time (both to perform testing within the time range of the original dataset and forecasting beyond its upper bound) is performed following the spectral decomposition of $\tilde{\mathbb{A}}$, i.e.
\begin{equation}
    \tilde{\mathbb{A}} \mathbb{W} = \mathbb{W}\Lambda
\end{equation}
with $\mathbb{W}$ being the eigenvectors matrix and the diagonal matrix $\Lambda$ containing the associated eigenvalues, which embed the dynamical evolution of the system. The former are needed to compute the so-called high-dimensional DMD modes $\boldsymbol{\Phi}\in\mathbb{R}^{\mathcal{N}_h\times r}$, i.e.
\begin{equation}
    \boldsymbol{\Phi} = \mathbb{X}^+\widetilde{\mathbb{V}}\widetilde{\Sigma}^{-1}\mathbb{W}
\end{equation}
which can be proven to be the eigenvectors of $\mathbb{A}$ \cite{Tu2014}. 

All the previous operations are performed only once during the training (offline) phase; during the test (online) phase, the solution can be reconstructed for unknown time values. Through expansion, the spectral decomposition in Eq. 6 allows for retrieving a simple expression for the advancement in time:
\begin{equation}
\label{advancement_in_time}
    \mathbf{x}_k \approx \sum_{j=1}^r\boldsymbol{\phi}_j \lambda_j^{k-1}b_j
\end{equation}
given $\boldsymbol{\phi}_j\in\mathbb{R}^{\mathcal{N}_h}$ the $j$-th column of $\boldsymbol{\Phi}$, $\lambda_j$ the $j$-th diagonal element of $\Lambda$ and $b_j$ the mode amplitude, computed from the projection of the initial condition onto the reduced space, i.e. $\mathbf{b}= \boldsymbol{\Phi}^T\mathbf{x}_1\in\mathbb{R}^r$. 

\subsection{Optimised version}\label{sec: opt-dmd}

One of the most significant limitations of the basic DMD approach is its susceptibility to noise, which can distort the relationship between consecutive state vectors and thus cause a wrong derivation of the spatio-temporal modes. There are several variants of the DMD algorithm that address this weakness, but over the years, the Optimised DMD \cite{optdmd2018} has been proven to be the most stable, unbiased and robust regression to fit the data \cite{bopdmd_2022}.  As explained in the previous section, the DMD aims at approximating the dynamics through a linear system, i.e. governed by laws of the kind $\der{\vec{x}}{t}\approx \mathbb{A}\vec{x}$, whose exact solution can be written as $\vec{x}=  e^{\mathbb{A}t}\,\vec{x}_0$. The evolution strongly depends on the eigenvalue-eigenvector pairs $\boldsymbol{\omega}-\boldsymbol{\Phi}$ of the dynamical matrix $\mathbb{A}$. Therefore, the general solution can be expressed as $\vec{x}=  \boldsymbol{\Phi} e^{\boldsymbol{\omega} t}\,\vec{x}_0$. With this in mind, the optimised version of DMD solves the following optimisation problem to perform a regression to exponential-time dynamics:
\begin{equation}
    \argmin\limits_{\boldsymbol{\omega}, \boldsymbol{\Phi}, \vec{b}} \norma{\mathbb{X} - \sum_{j=1}^rb_j\boldsymbol{\phi}_j e^{\omega_j t}} =  \argmin\limits_{\boldsymbol{\omega}, \boldsymbol{\Phi}, \vec{b}} \norma{\mathbb{X} - \boldsymbol{\Phi}\text{exp}(\boldsymbol{\omega}t)\vec{b}}
    \label{eqn: opt-dmd-pb}
\end{equation}
where an $r$-rank exponential approximation of the snapshots is performed. In Eq. \eqref{eqn: opt-dmd-pb}, $\boldsymbol{\Phi}$ and $\boldsymbol{\phi}_j$ are the DMD modes, $\omega_j$ are the continuous eigenvalues associated to the discrete eigenvalues $\boldsymbol{\lambda}$ \cite{brunton_data-driven_2022} and $\vec{b}$ represents the amplitudes. This expansion follows the spectral decomposition of the operator $\mathbb{A}$ discussed in the previous section and carries a significant advantage compared to the basic DMD algorithm. Whereas the latter requires the sampling of spatio-temporal data to be uniform in time, the optimised version of DMD allows the generalisation of this process for non-uniform trajectories. This reasoning resulted in the Bagging version of the Optimised DMD \cite{bopdmd_2022} to treat noisy data from experiments, where subsets of the original datasets are used to create an ensemble of optimised DMD models, whose outputs are then averaged.

\section{Parametric Dynamic Mode Decomposition}\label{sec: pdmd}

In the previous section, the DMD algorithm for a single-parameter snapshot matrix has been presented. However, this algorithm lacks the capability of dealing with parametric datasets, which is a desired feature for its application to engineering systems. To the best of the authors' knowledge, three different versions of parametric DMD exist in the literature:
\begin{itemize}
    \item Reduced Operators Interpolation (ROI): this version has been proposed by the authors of the present work in \cite{NUTHOS24_pDMD}, and it is based on the search for common structures in the (parametric) dynamical matrices themselves (Section \ref{sec: pdmd-nuthos}).
    \item Reduced Koopman Operator Interpolation (RKOI): this version derives from the original work of Schmid \cite{schmid_dynamic_2010}, which has been discussed in \cite{huhn_parametric_2023} and extended to the optimized version in \cite{faraji_data-driven_2024-2}; it is based on the interpolation of the eigenvalues/eigenvectors pairs of single-parameter DMD models (Section \ref{sec: pdmd-kutz}).
    \item Interpolation of the Latent Space: this version, proposed in \cite{andreuzzi_dynamic_2023}, is part of the pyDMD package \cite{demo_pydmd_2018, ichinaga_pydmd_2024}; it is based on the interpolation of the latent dynamics, according to two different strategies (Section \ref{sec: pdmd-rozza}), called Monolithic (Mono) and Partitioned (Part).  
\end{itemize}

For all versions, different snapshot matrices $\mathbb{X}^{\boldsymbol{\mu_i}} = [\mathbf{x}^{\boldsymbol{\mu_i}}(t_1) \mid \mathbf{x}^{\boldsymbol{\mu_i}}(t_2) \mid \dots \mid \mathbf{x}^{\boldsymbol{\mu_i}}(t_{N_t})] \in \mathbb{R}^{\mathcal{N}_h \times N_t}$ are collected for each parameter $\boldsymbol{\mu}_i$ within the parametric domain $\mathcal{D}\in\mathbb{R}^p$ ($p\geq 1$). If the SVD is performed separately on each matrix, a set of modes for each parameter is obtained, making them parameter-dependent: thus, a suitable interpolant for them must be implemented within the algorithm, as done in \cite{huhn_parametric_2023}, however, this can be drastically inefficient if high-dimensional modes are to be learnt. Instead, the snapshots can be stacked in a single matrix $\mathcal{X}\in\mathbb{R}^{\mathcal{N}_h \times N_t\cdot N_p}$, defined as follows (with $N_p$ the number of parameters), i.e.
\begin{equation}
    \mathcal{X} = [\mathbb{X}^{\boldsymbol{\mu_1}} \mid \mathbb{X}^{\boldsymbol{\mu_2}} \mid \dots \mid \mathbb{X}^{\boldsymbol{\mu_{N_p}}}]
\end{equation}
and the randomised version of the SVD, for computational efficiency \cite{halko2010findingstructurerandomnessprobabilistic}, can be adopted to find a set of modes $\widetilde{\mathbb{U}}\in\mathbb{R}^{\mathcal{N}_h \times r}$, which encodes the most dominant spatial structures of the whole dataset. In this way, each parametric matrix is projected onto a reduced space:
\begin{equation}
    \widetilde{\mathbb{V}}^{\boldsymbol{\mu_i}} = \widetilde{\mathbb{U}}^T\mathbb{X}^{\boldsymbol{\mu_i}}
\end{equation}
with $\widetilde{\mathbb{V}}^{\boldsymbol{\mu_i}}=[\mathbf{v}_1^{\boldsymbol{\mu_i}}, \mathbf{v}_2^{\boldsymbol{\mu_i}}, \dots, \mathbf{v}_{N_t}^{\boldsymbol{\mu_i}}]\in\mathbb{R}^{r\times N_t}$.

According to this procedure, the DMD algorithms are applied to the latent dynamics in the compressed space rather than working at the high-dimensional level. This approach can drastically reduce the computational time needed to generate the DMD operators, even though some numerical instabilities may arise \cite{andreuzzi_dynamic_2023} when the number of rows (dimension of the state vector, either the spatial degrees of freedom or the rank of the SVD) is much lower than the number of time instances: a possible remedy for this problem is provided by higher-order variants of DMD \cite{LeClainche_HoDMD}. As already highlighted in \cite{faraji_data-driven_2024-2}, in parametric DMD, the choice of the optimal rank $r$ of the SVD involves a compromise between two aspects: the information content captured by the SVD modes in the latent representation and the accuracy of the interpolation in the parameter space on the test data points. In particular, increasing the number of SVD modes improved the closeness of the predictions of a single DMD to the original training data corresponding to each specific training parameter. However, the interpolation quality between the DMD bases degrades as the number of DMD modes increases: the optimal value of $r$ depends on the specific application and on its main objective (that is, interpolation, reconstruction or extrapolation).

The ROI and RKOI algorithms have been implemented by the authors, and the associated classes are available on \href{https://github.com/ERMETE-Lab/pDMD}{Github} (the generation of the single-parameter DMD models is however made through the pyDMD package \cite{demo_pydmd_2018, ichinaga_pydmd_2024}); the version of \cite{andreuzzi_dynamic_2023} can be found in the pyDMD package itself. In the next Sections, the three algorithms will be explained in more detail; summarising schemes of each method can be found in Appendix \ref{app: schemes-pdmd}.

\subsection{Reduced Operators Interpolation}\label{sec: pdmd-nuthos}

This approach, referred to hereafter as Reduced Operators Interpolation (ROI), has been proposed by the authors in \cite{NUTHOS24_pDMD}, assuming that the DMD operators may share some dominant structures among different parameter realisations. Therefore, for each parametric snapshot matrix $\mathbb{X}^{\boldsymbol{\mu}_i}$, a DMD model, adopting the basic DMD algorithm (Section \ref{sec: basic-dmd}), is generated from the SVD coefficients such that the dynamical system in the latent space evolves as 
\begin{equation}
    \mathbf{v}_{k+1}^{\boldsymbol{\mu_i}} = \tilde{\mathbb{A}}^{\boldsymbol{\mu_i}} \mathbf{v}_{k}^{\boldsymbol{\mu_i}}
    \label{eqn: roi-time-dynamics}
\end{equation}
In this way, a collection of operators $\tilde{\mathbb{A}}^{\boldsymbol{\mu_i}}\in\mathbb{R}^{r\times r}$ can be generated. These operators are unfolded into column vectors and collected in a new 'snapshot' matrix $\mathcal{A}\in\mathbb{R}^{r^2\times N_p}$ defined as
\begin{equation}
    \mathcal{A} = [\tilde{\mathbf{a}}^{\boldsymbol{\mu_1}} \mid \tilde{\mathbf{a}}^{\boldsymbol{\mu_2}} \mid \dots \mid \tilde{\mathbf{a}}^{\boldsymbol{\mu_{N_p}}}]
\end{equation}
being $\tilde{\mathbf{a}}^{\boldsymbol{\mu_i}}$ the $i$-th unfolded operator correspondent to $\tilde{\mathbb{A}}^{\boldsymbol{\mu_i}}$. Matrix $\mathcal{A}$ then collects the coefficients of each underlying linear operator. It is then legitimate to investigate if there is some dominant behaviour among them for different parametric instances and if there is some relationship describing their behaviour as a function of the parameters. To discover this relationship, a second SVD is performed on $\mathcal{A}$ so that each operator within the parametric range can be expressed as:
\begin{equation}
    \tilde{\mathbb{A}}^{\boldsymbol{\mu_i}} \simeq \sum_{k=1}^{r_a}\xi_k(\boldsymbol{\mu_i})\cdot \mathbb{F}_k
    \label{eqn: svd-op}
\end{equation}
Hence, a set of SVD modes for the dynamical matrix is encoded in $\{\mathbb{F}_k\}_{k=1}^{r_a}$, whereas the parametric variation is embedded in the coefficients $\{\xi_k(\boldsymbol{\mu_i})\}_{k=1}^{r_a}$, given $r_a$ the rank of this second SVD\footnote{According to experience, it is suggested to keep this value close to the number of parameters evaluation $N_p$ to ensure a good reconstruction of the operator, since the DMD algorithm is very sensitive to it \cite{NUTHOS24_pDMD}.}. In this way, the dependence on the parameters is detached from the operators; the last step of the algorithm is to learn this parametric dependence of the operators using a regression model aptly chosen. The regressor $\mathcal{I}$ is trained on the reduced coefficients $\xi_k(\boldsymbol{\mu})$ of the operators and the training parameter values $\boldsymbol{\mu}$, thus obtaining a surrogate map $\mathcal{F}_k$ that can predict the operators for unseen parameter values:
\begin{equation}
    \boldsymbol{\mu}\xrightarrow{\mathcal{F}_k} \xi_k(\boldsymbol{\mu})
\end{equation}
Adopting this approach, in principle, any operator within the parametric range can be inferred. The use of a second SVD significantly reduces the number of interpolants that need to be computed; conversely, the method developed by Huhn et al. \cite{huhn_parametric_2023} directly interpolates each element of the DMD operators, meaning that for high values of $r$, a large number of interpolants would be required. This problem may also affect the RKOI method, discussed in Section \ref{sec: pdmd-kutz}.

Once the training phase is completed, a new parametric transient can be simulated for $\boldsymbol{\mu}^\star$: the dynamical matrix $\tilde{\mathbb{A}}^{\boldsymbol{\mu}^\star}$ can be estimated by evaluating the regression model $\{\mathcal{F}_k\}_{k=1}^{r_a}$ at $\boldsymbol{\mu}^\star$ to retrieve the coefficients $\xi_k(\boldsymbol{\mu}^\star)$, from which the unknown matrix can be approximated using Eq. \eqref{eqn: svd-op}. The system is then advanced in time using Eq. \eqref{eqn: roi-time-dynamics}, from which the latent dynamics $\widetilde{\mathbb{V}}^{\boldsymbol{\mu}^\star}$ of the new parameters are calculated. By projecting back onto the reduced basis spanned by $\widetilde{\mathbb{U}}$, i.e. $\mathbb{X}^{\boldsymbol{\mu}^\star}\simeq \widetilde{\mathbb{U}}\widetilde{\mathbb{V}}^{\boldsymbol{\mu}^\star}$, the overall high-dimensional state vector can be finally retrieved.

\subsection{Reduced Koopman Operator Interpolation}\label{sec: pdmd-kutz}

In the same way as the ROI version, the Reduced Koopman Operator Interpolation (RKOI) version of parametric DMD \cite{faraji_data-driven_2024-2, huhn_parametric_2023} starts by building a DMD model of the latent dynamics with the Optimised DMD \cite{optdmd2018} for each parametric configuration, obtaining (expanding in the continuous time domain $t$):
\begin{equation}
    \mathbf{v}^{\boldsymbol{\mu_i}}(t)= \sum_{j=1}^r\boldsymbol{\phi}_j^{\boldsymbol{\mu_i}} \cdot  \text{exp}(\omega_j^{\boldsymbol{\mu_i}}\cdot t)\cdot  b_j^{\boldsymbol{\mu_i}}
    \label{eqn: dmd-decomp-continuous}
\end{equation}
with $\boldsymbol{\phi}^{\boldsymbol{\mu_i}}_j\in\mathbb{R}^r$ being the $j$-th DMD mode, $\omega_j^{\boldsymbol{\mu_i}}$ being the $j$-th eigenfrequency, descending from the eigenvalues of $\tilde{\mathbb{A}}^{\boldsymbol{\mu_i}}$, and $b_j^{\boldsymbol{\mu_i}}$ being the mode amplitude computed from the projection of the initial condition onto the reduced space, i.e. $\mathbf{b}^{\boldsymbol{\mu_i}}= \boldsymbol{\Phi}^{*,\boldsymbol{\mu_i}}\mathbf{v}^{\boldsymbol{\mu_i}}(t_0)$, the superscript $^*$ indicating the conjugate transpose.

Each element of Eq. \eqref{eqn: dmd-decomp-continuous} is dependent on the parameter $\boldsymbol{\mu}$. The general representation of  $\mathbf{v}(\boldsymbol{\mu})$ as a function of the parameter can be retrieved by building regression models for the DMD modes, the frequencies and the amplitudes. Each regression model is trained on the reduced DMD modes $\boldsymbol{\phi}_j^{\boldsymbol{\mu_i}}$, the frequencies $\omega_j^{\boldsymbol{\mu_i}}$ and the amplitudes $b_j^{\boldsymbol{\mu_i}}$ from the training dataset and the correspondent parameter values $\boldsymbol{\mu}$, thus obtaining a map that can predict the DMD modes, frequencies and amplitudes for unseen parameter values. As highlighted in \cite{faraji_data-driven_2024-2}, it is remarked that this interpolation would result in a reasonable estimation of the dynamics only if the variation of the spatio-temporal modes is smooth across the parameter space. If there are strong non-linear dependencies, they can cause significant discrepancies between the interpolated and the 'true' spatial and temporal characteristics.

Once the latent dynamics for the unseen parameter $\boldsymbol{\mu}^\star$ are retrieved, the high-dimensional state vector can then be reconstructed using the spatial SVD modes:
\begin{equation}
    \mathbf{x}^{\boldsymbol{\mu^\star}}(t) \simeq \widetilde{\mathbb{U}}\mathbf{v}^{\boldsymbol{\mu^\star}}(t) =  \sum_{j=1}^r \left(\sum_{i=1}^r \mathbf{u}_i\cdot {\phi}_{j,i}^{\boldsymbol{\mu^\star}} \right) \cdot  \text{exp}(\omega_j^{\boldsymbol{\mu_i}}\cdot t) \cdot b_j^{\boldsymbol{\mu^\star}}
\end{equation}
being $\mathbf{u}_i\in\mathbb{R}^{\mathcal{N}_h}$ the $i$-th SVD mode, that is, the $i$-th column of $\widetilde{\mathbb{U}}$. This version of parametric DMD is quite similar to ROI; nevertheless, the use of the Optimised DMD algorithm allows for better catching the temporal dynamics, thus reducing the dimension of the latent space. This version, presented in \cite{faraji_data-driven_2024-2}, can be seen as an extension of \cite{huhn_parametric_2023}.

\subsection{Interpolation of the Latent Space}\label{sec: pdmd-rozza}

As in the two previous versions, the starting point of this latter method for parametric DMD consists in applying a SVD-based reduction to the parameter-dependent snapshots to obtain their reduced representation; then, non-parametric DMD is applied to each snapshot realisation to approximate their evolution for future time instants. The parametric variability is then explored using regression techniques in the latent space. Two approaches for applying DMD to the reduced snapshots are described in \cite{andreuzzi_dynamic_2023}: in the first one, the \textbf{monolithic} approach, a unique linear operator is generated to fit the dynamics of the entire (parametrized) system; alternatively, in the \textbf{partitioned} approach $N_p$ linear operators are constructed to approximate the dynamics of the parametric configurations.

\subsubsection{Monolithic Approach}

In this approach, the reduced dynamics are stored in the matrices $\widetilde{\mathbb{V}}^{\boldsymbol{\mu_i}}$ for each parameters $\boldsymbol{\mu}_i$; these matrices are stacked together into a matrix $\mathcal{V}_2\in\mathbb{R}^{r\cdot p\times N_t}$: 
\begin{equation}
    \mathcal{V}_2 = \left[
        \begin{array}{c}
             \widetilde{\mathbb{V}}^{\boldsymbol{\mu_1}}\\
             \widetilde{\mathbb{V}}^{\boldsymbol{\mu_2}}\\
             \vdots \\
             \widetilde{\mathbb{V}}^{\boldsymbol{\mu_p}}
        \end{array}
    \right]
\end{equation}
whose columns are the temporal trajectories. An overall DMD operator $\widetilde{A}\in\mathbb{R}^{r\cdot p\times r\cdot p}$ can be obtained to advance in time the overall stacked reduced state vector:
\begin{equation}
    \left[
        \begin{array}{c}
             \vec{v}_{k+1}^{\boldsymbol{\mu_1}}\\
             \vec{v}_{k+1}^{\boldsymbol{\mu_2}}\\
             \vdots\\
             \vec{v}_{k+1}^{\boldsymbol{\mu_p}}
        \end{array}
    \right] = \widetilde{A}\left[
        \begin{array}{c}
             \vec{v}_{k}^{\boldsymbol{\mu_1}}\\
             \vec{v}_{k}^{\boldsymbol{\mu_2}}\\
             \vdots\\
             \vec{v}_{k}^{\boldsymbol{\mu_p}}
        \end{array}
    \right] 
    \label{eqn: monolithic-op}
\end{equation}

Adopting this approach, the dynamics of the parametric system are expressed by a single DMD operator, allowing for the detection of recurrent patterns in the dynamics of different parametric configurations \cite{andreuzzi_dynamic_2023}. All these operations are to be performed during the offline phase; in the online phase, this version of parametric DMD adopts an approach very similar to that of POD with Interpolation (POD-I) \cite{demo_complete_2019, ortali_gaussian_2022}, in which the SVD space is approximated with some regression or interpolation techniques. Let $\mathcal{I}:\mathcal{D}\rightarrow \mathbb{R}^r$ be the regressor taking the parameter $\boldsymbol{\mu}$ as input and returning the reduced coefficients at the desired time $t$. This regressor is trained using all the reduced dynamics obtained by solving the trained DMD monolithic model. To calculate the approximated reduced snapshot, a multidimensional regressor is needed; alternatively, a single regressor can be built for single components of the reduced snapshots. Once the latent dynamics have been calculated, the high-dimensional state can be retrieved by projecting back onto the SVD modes. Thus, the training of the regressor occurs during the online phase, as it requires the solution of the DMD models: from the computational point of view, this approach will be more expensive compared to the other versions of parametric DMD.

\subsubsection{Partitioned Approach}

This last approach combines the ROI version and the previous monolithic procedure. Instead of building a monolithic DMD operator, $N_p$ separate DMDs, one for each matrix $\widetilde{\mathbb{V}}^{\boldsymbol{\mu_i}}$, are performed as in the ROI approach. The main differences are in the online phase: in fact, for a new time instant $t^*$, the DMD models are solved independently, and the reduced coefficients are interpolated with a regressor $\mathcal{I}$ as in the monolithic approach.

\section{Numerical Results}\label{sec: num-res}

The three different versions of the parametric DMD have been applied to three different test cases: a flow over cylinder for low Reynolds number solved with \textit{dolfinx} (v. 0.6.0) \cite{ alnaes_unified_2014, baratta_dolfinx_2023, scroggs_basix_2022, scroggs_construction_2022} within the OFELIA framework developed by the authors in \cite{loi_ofelia_2024}; a benchmark dataset for the flow over cylinder problem, designed for machine learning testing with CFD data \cite{luo_cfdbench_2024}; the DYNASTY experimental facility \cite{cammi_dynasty_2016, benzoni_preliminary_2023} deployed at Politecnico di Milano, meant for investigating natural circulation phenomena of internally heated fluids in the framework of Generation-IV reactors \cite{generation_iv_international_forum_technology_2014}.

The case studies selected in this paper are aimed at verifying the accuracy, the advantages and the shortcomings of the different pDMD algorithms available in the literature, starting at first from rather simple fluid dynamics problems yet presenting complex dynamics and significant non-linearities which are generally not trivial to predict (i.e., flow over cylinder); then, the case study using the data R5 model for the DYNASTY experimental facility has been designed to show that these methods can also work on real system, using as a starting model a system code widely adopted in the nuclear community.

\subsection{Laminar Flow Over Cylinder}

\begin{figure}[htbp]
    \centering
    \includegraphics[width=0.75\linewidth]{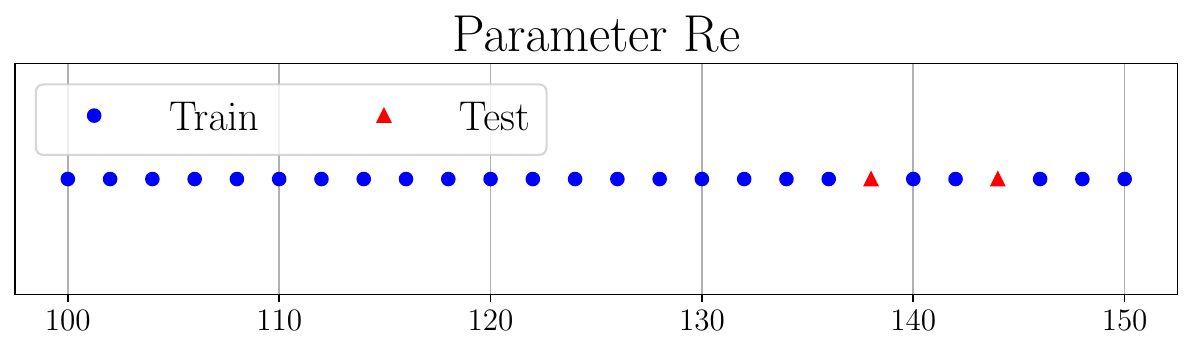}
    \caption{Train and Test split of the Reynolds number for the Flow Over Cylinder solved with \textit{dolfinx}.}
    \label{fig: flow-fenics-traintest}
\end{figure}

The 'flow over cylinder' is one of the most adopted test cases for thermal-hydraulics, ROM and ML techniques \cite{faraji_data-driven_2024-1, faraji_data-driven_2024-2, andreuzzi_dynamic_2023, luo_cfdbench_2024, lorenzi_pod-galerkin_2016, stabile_pod-galerkin_2017}, as even for low values of the Reynolds number (i.e., in the laminar flow regime), the dynamics are non-trivial due to the oscillations of the vortex shedding. For the first test case, a channel of width 1.1 m with an obstacle of radius $\delta=0.05$ m placed at $y=0.5$ m is considered (from the lower wall). More details about the test case and the dataset can be found on the GitHub repository associated with this work (see Code and Data Availability). This problem is governed by the incompressible Navier-Stokes equations
\begin{equation}
    \left\{
        \begin{aligned}
            &\dpart{\vec{u}}{t}+(\vec{u}\cdot\nabla)\vec{u}=\nu\Delta \vec{u}-\nabla p\\
            &\nabla\cdot \vec{u}=0
        \end{aligned}
    \right.
\end{equation}
with $\vec{u}$ and $p$ being, respectively, the velocity and pressure and $\nu$ being the kinematic viscosity. The problem is closed by suitable boundary conditions: constant parabolic velocity at the inlet, directed along $x$, null pressure at the outlet and no-slip conditions at the walls. This non-linear system of PDEs has been discretised in space using the Finite Element method with the \textit{dolfinx} \cite{baratta_dolfinx_2023, scroggs_basix_2022, scroggs_construction_2022, alnaes_unified_2014} library for Python, using the code from \cite{loi_ofelia_2024}. The parameter $\boldsymbol{\mu}$ for this problem is the Reynolds number, defined as $Re = \frac{1}{\nu}\cdot (u_{x, in}\cdot 2\delta)$: 26 simulations have been performed in the range $Re \in [100,150]$, considering equally spaced parameter values. Given the sensitivity of parametric techniques to the completeness of the training set, the parameter space has been split manually, taking only two parameter realisations at high values of $Re$ (where the flow is expected to be more complex) for the test set; this subdivision is reported in Figure \ref{fig: flow-fenics-traintest}.

The snapshots matrix $\mathcal{X}$ for this problem contains the velocity field $\vec{u}$: it has dimensions $N_p = 26 \times \mathcal{N}_h = 4951 \times N_t = 1000$, respectively the number of parameters, the spatial degrees of freedom (including both $x$ and $y$ components) and the number of equally spaced time instances. Regardless of the parametric DMD version, the first step consists of performing the SVD on the training snapshots to obtain the SVD modes encoding the dominant spatial structures and the singular values weighting the information retained by each mode: the randomised version \cite{halko2010findingstructurerandomnessprobabilistic} of the SVD algorithm is used for this purpose.

\begin{figure}[htbp]
    \centering
    \includegraphics[width=1\linewidth]{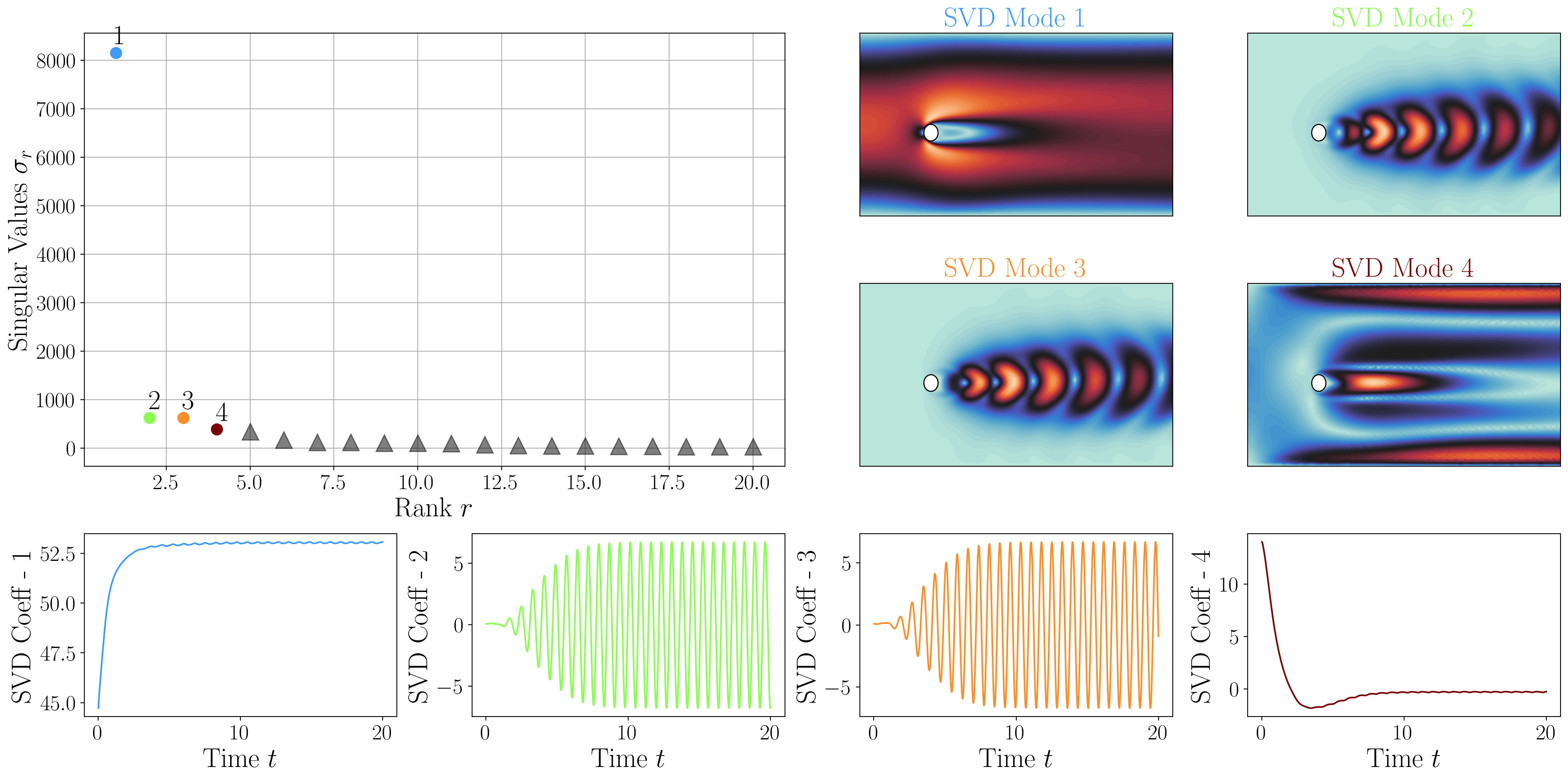}
    \caption{Decay of the singular values (top left), contour plots of the first 4 SVD modes of velocity $\mathbf{u}$ (top right), underlining the hierarchical spatial features and temporal evolution of the latent dynamics for fixed Reynolds (bottom). The singular values show an exponential decay, highlighting the fact that low rank modes retain the majority of the information content.}
    \label{fig: flow-fenics-modes}
\end{figure}

Figure \ref{fig: flow-fenics-modes} shows the decay of the singular values (on the top left), the first 4 SVD spatial modes for the velocity field (on the top right) and the dynamics of the latent space (on the bottom): the first mode encodes the average flow field, whereas from the second one onward the low-scale dynamics are embedded; in particular, the vortex shedding oscillations can be mainly described by the 2nd and 3rd spatial modes, whereas higher-order modes contain higher-order phenomena, as underlined by the associated SVD coefficients. 

The different algorithms for parametric DMD will be compared in terms of computational time and average relative error measured with the Frobenius norm: let $\mathbb{X}^{\boldsymbol{\mu}^\star}$ be the snapshot matrix for the test parameter $\boldsymbol{\mu}^\star$ and let $\hat{\mathbb{X}}_*^{\boldsymbol{\mu}^\star}$ be the DMD reconstruction (for a specific version $*$ of the pDMD); then, the error $\epsilon^{\boldsymbol{\mu}^\star}_*$ can be calculated as:
\begin{equation}
    \epsilon^{\boldsymbol{\mu}^\star}_* = \frac{\norma{\mathbb{X}^{\boldsymbol{\mu}^\star}-\hat{\mathbb{X}}_*^{\boldsymbol{\mu}^\star}}_F}{\norma{\mathbb{X}^{\boldsymbol{\mu}^\star}}_F}
    \label{eqn: rec-error}
\end{equation}

To make a fair comparison in terms of computational times and errors, at first, the rank $r=10$ is chosen for all the pDMD versions. This value can be selected by looking at the decay in the singular values of Figure \ref{fig: flow-fenics-modes}, following the fact that each singular value describes the relative information retained by the modes \cite{lassila_model_2014, rozza_model_2020, brunton_data-driven_2022}: at rank 10, the exponential decay has already been exhausted, and no meaningful information is added by retain a higher number of modes\footnote{As already mentioned, the rank $r$ is an hyperparameter of DMD, and its optimal value depends on the application and on the main objective of the analysis. Choosing a value of $r$ for which, visually, the decay of the singular values reaches a plateau is one way to provide an initial guess for it.} Additionally, all algorithms have been run on the same hardware with Intel Core i7-9800X CPU and a clock speed of 3.80 GHz.

\begin{table}[htbp]
    \centering
    \begin{tabular}{c|c|c|c}
        \hline
        Parameter $\boldsymbol{\mu}^\star$ & Algorithm $*$ & Error $\epsilon^{\boldsymbol{\mu}^\star}_*$ & CPU Time (s) \\ 
        \hline
        \multirow{4}{*}{138.0} 
        & ROI & \cellcolor{red!20}0.1548 & \cellcolor{green!20}0.2227 \\ 
        & RKOI & 0.0648 & 1.3001 \\ 
        & Monolithic & \cellcolor{green!20}0.0390 & 15.3591 \\ 
        & Partitioned & 0.1541 & \cellcolor{red!20}17.9683 \\ 
        \hline
        \multirow{4}{*}{144.0} 
        & ROI & \cellcolor{red!20}0.1580 & \cellcolor{green!20}0.6157 \\ 
        & RKOI & 0.0679 & 1.2875 \\ 
        & Monolithic & \cellcolor{green!20}0.0421 & \cellcolor{red!20}18.3177 \\ 
        & Partitioned & 0.1579 & 14.8326 \\ 
        \hline
    \end{tabular}
    \caption{Error (as defined in Eq. \eqref{eqn: rec-error}) and CPU Time for different pDMD algorithms for the two parameter values of the test set and same value of the rank.}
    \label{tab:flow-fenics-same-rank-comparison}
\end{table}

The different pDMD versions are compared in terms of errors (computed according to Eq. \eqref{eqn: rec-error}) and CPU time (in seconds) in Table \ref{tab:flow-fenics-same-rank-comparison}. The ROI and Partitioned approaches show the worst performance in terms of errors, as the adopted rank is not enough to correctly build the reduced model: whereas the mean flow is well described by these two approaches, the vortices are not correctly predicted (see the Supplementary materials for the videos of the temporal evolution of the reconstructed velocity field using the various pDMD versions and its residual). On the other hand, the RKOI and the Monolithic approaches work better: the former is built using the optimised version of the DMD, which is better for learning the dynamics and more stable when the size of the state (i.e., the rank) is much lower than the number of time steps; the latter, instead, works with a higher size of the state vector since a single operator is built on the stacked training latent dynamics, as in Eq. \eqref{eqn: monolithic-op}. In terms of computational time, however, both the versions implemented in the pyDMD package (Monolithic and Partitioned) are more expensive since they fit the regressors within the online phase; on the contrary, both ROI and RKOI move this task in the offline stage, thus allowing for better computational performances\footnote{Regarding CPU times, the computational bottleneck is given by the interpolant step: flows at higher Reynolds number, being more complex in terms of dynamics, require more computational time as the interpolant is more complex}.

After this comparison at a fixed value of the rank, each version of the pDMD has been trained by tuning the rank of the latent dynamics to have, for all versions, a comparable error in the training test. The ROI requires a rank value of $40$ to obtain accurate results on the training data; for RKOI, a rank of 10 is sufficient due to its better capabilities in building the surrogate model at a lower rank; for the Monolithic approach, a rank of 8 for each parameter is sufficient (note that the state vector in the latent representation would be $8\cdot 24$); finally, the Partitioned approach has been trained with rank equal to 30. Except for the RKOI and Monolithic approaches, "higher" ranks have been used to push the algorithms towards optimal performance (i.e., as-best-as-possible): as shown before, if these were lowered, the accuracy of the DMD prediction, even for the training parameters, would worsen a lot.

\begin{table}[htbp]
    \centering
    \begin{tabular}{c|c|c|c}
        \hline
        Parameter $\boldsymbol{\mu}^\star$ & Algorithm $*$ & Error $\epsilon^{\boldsymbol{\mu}^\star}_*$ & CPU Time (s) \\ 
        \hline
        \multirow{4}{*}{138.0} 
        & ROI & 0.0336 & \cellcolor{green!20}0.677 \\ 
        & RKOI & \cellcolor{red!20}0.0648 & 1.403 \\ 
        & Monolithic & 0.0462 & 17.867 \\ 
        & Partitioned & \cellcolor{green!20}0.0305 & \cellcolor{red!20}34.755 \\ 
        \hline
        \multirow{4}{*}{144.0} 
        & ROI & \cellcolor{green!20}0.0303 & \cellcolor{green!20}0.408 \\ 
        & RKOI & \cellcolor{red!20}0.0679 & 1.361 \\ 
        & Monolithic & 0.0526 & 11.060 \\ 
        & Partitioned & 0.0320 & \cellcolor{red!20}29.834 \\ 
        \hline
    \end{tabular}
    \caption{Relative error in the Frobenius norm and CPU times for different pDMD algorithms for the two parameter values of the test set, using optimal ranks for the different versions.}
    \label{tab: flow-fenics-residuals-time-opt-rank}
\end{table}

\begin{figure}[htbp]
    \centering
    \includegraphics[width=1\linewidth]{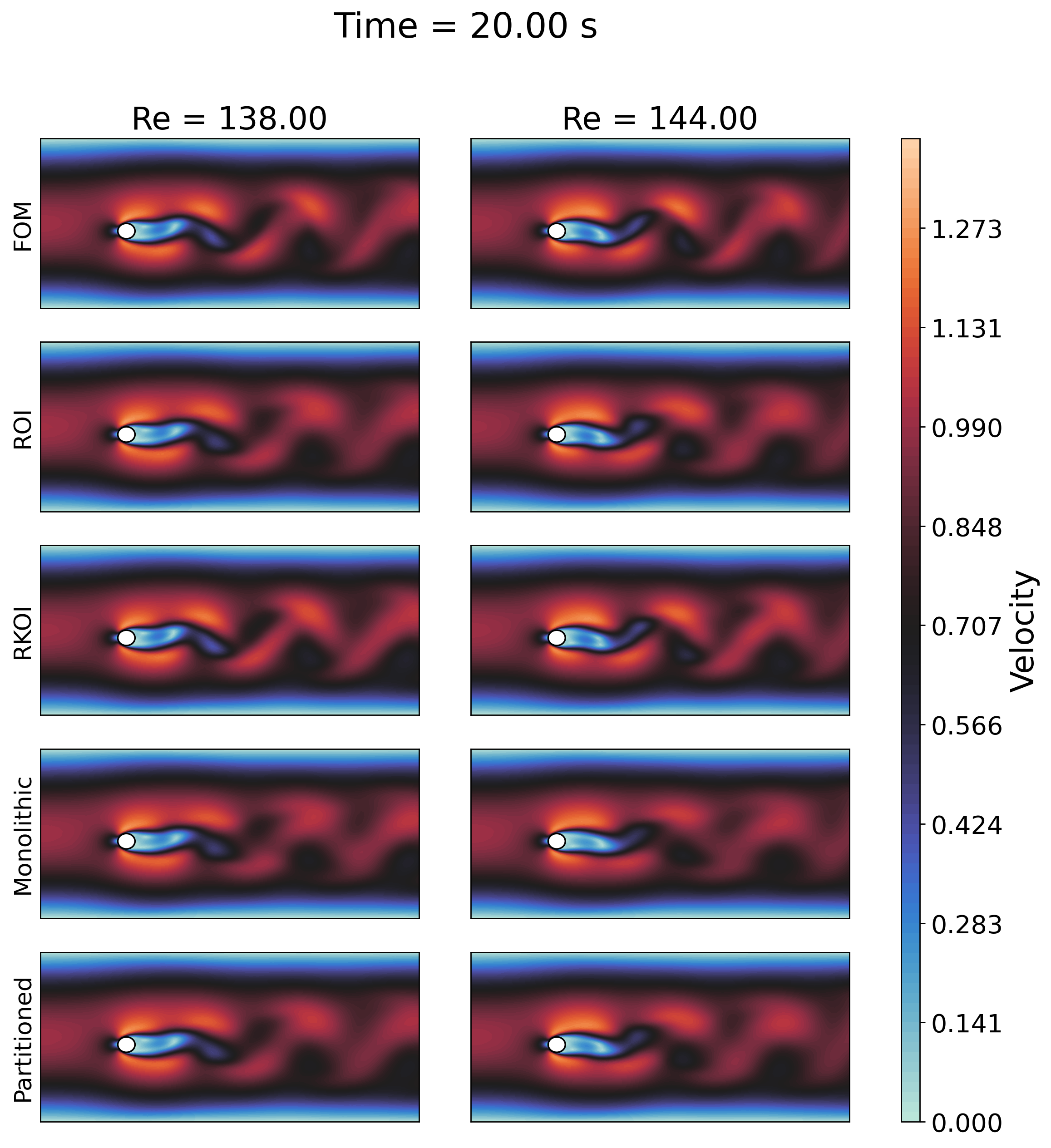}
    \caption{Contour plots of the velocity magnitude for the test parameters at the final time of the simulation, comparing the pDMD versions with the full-order solution (FOM).}
    \label{fig: flow-fenics-contour}
\end{figure}

After tuning the rank for each version, their performances are compared in terms of the error computed using Eq. \eqref{eqn: rec-error}, and the latter is reported in Table \ref{tab: flow-fenics-residuals-time-opt-rank}: if properly tuned in terms of the value of the rank, each method provides an accurate prediction of the flow field for unseen Reynolds numbers within the training range. The RKOI shows a slightly higher error, as it is less able to predict the behaviour of the velocity for the initial time instants, but it is better in predicting the oscillations, as shown in Figure \ref{fig: flow-fenics-contour}, in which the prediction of each method is plotted and compared with the FOM solution. The temporal evolution of the error for the three algorithms will be explored better in the next test case.

\subsection{Flow Over Cylinder (Benchmark Dataset)}

\begin{figure}[htbp]
    \centering
    \includegraphics[width=0.75\linewidth]{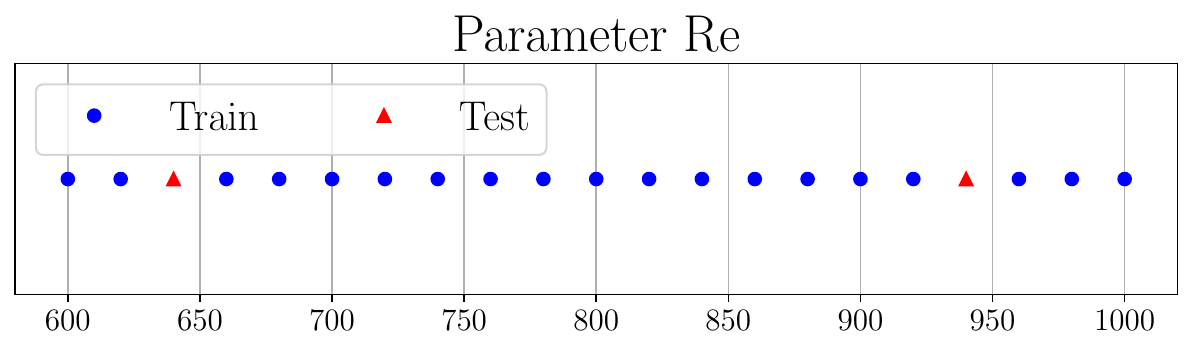}
    \caption{Train and Test split of the Reynolds number for the Flow Over Cylinder from the CFDbench dataset \cite{luo_cfdbench_2024}.}
    \label{fig: flow-cfdbench-traintest}
\end{figure}

A similar configuration of the previous case is now considered, using the CFDbench dataset \cite{luo_cfdbench_2024} dataset, which is purposefully designed to train and test machine learning algorithms. In this dataset, several case studies are solved for parametric configurations, both including physical and geometrical variations. The selected geometry is the same as in the previous case, albeit with different dimensions, and the flow physics are still governed by the Navier-Stokes equations. The inlet boundary condition for velocity now varies (uniformly in space) from 3 m/s to 5 m/s: given the thermo-physical properties of density $\rho = 10$ kg/m$^3$ and dynamic viscosity $\mu=1e-3$ Pa s, and an hydraulic diameter of 0.02 m, the parameter range in terms of the Reynolds number is now $Re\in[600,1000]$. For each parameter realisation, 1000 time steps have been saved and the data are stored as images shaped $64\times 64$ pixels. The parameters have been randomly divided into train and test as in Figure \ref{fig: flow-cfdbench-traintest}. 

Building from the results obtained previously, the different DMD versions will be compared, considering for each the optimal value of the rank. This test case serves as an investigation of the use of parametric DMD on snapshots coming from synthetic images; nevertheless, frames coming from experiment recordings can also be used (for instance, the frames can be the output in time of a thermo-camera or particle image velocimetry), removing the need of running numerical simulations and instead generating a surrogate model of the observed physics in a model-free manner. Regardless, the relative error over time $\varepsilon^{\boldsymbol{\mu}^\star}_*({t})$ is defined as:
\begin{equation}
    \varepsilon_*^{\boldsymbol{\mu}^\star}({t}) = \frac{\norma{\vec{x}^{\boldsymbol{\mu}^\star}(t)-\hat{\vec{x}}_*^{\boldsymbol{\mu}^\star}(t)}_2}{\norma{\vec{x}^{\boldsymbol{\mu}^\star}(t)}_2}
    \label{eqn: rec-time-error}
\end{equation}
where $\vec{x}^{\boldsymbol{\mu}^\star}(t)$ represents the high-dimensional state at time $t$ for the parameter $\boldsymbol{\mu}^\star$ with $\hat{\vec{x}}_*^{\boldsymbol{\mu}^\star}(t)$ as the DMD prediction for the different algorithms.

\begin{figure}[htbp]
    \centering
    \includegraphics[width=1\linewidth]{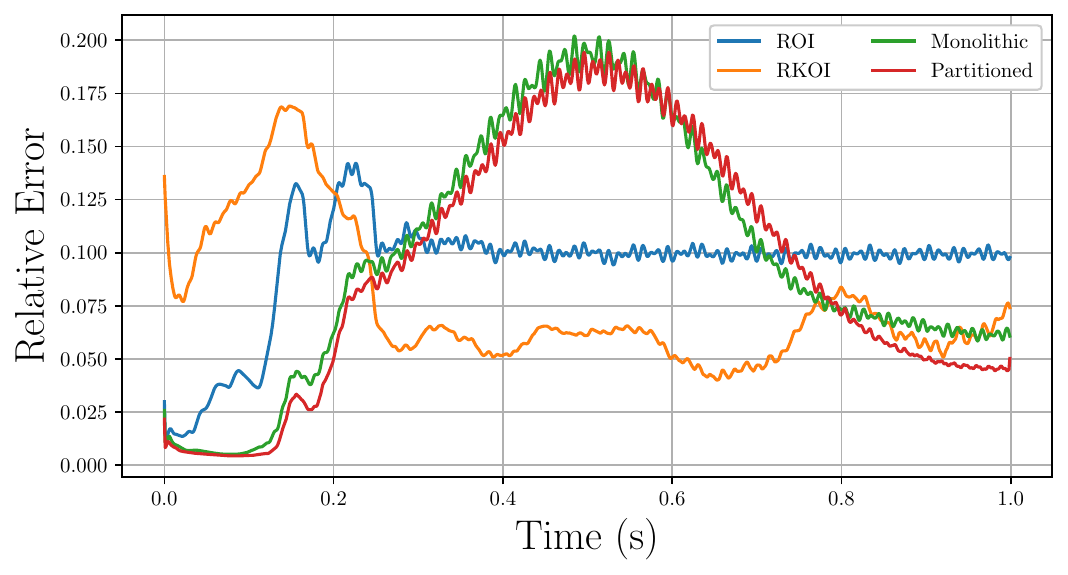}
    \caption{Average (with respect to test parameters) relative errors over time for the different pDMD versions applied on the CFDbench dataset.}
    \label{fig: flow-cfdbench-temporal-errors}
\end{figure}

\begin{figure}[htbp]
    \centering
    \includegraphics[width=1\linewidth]{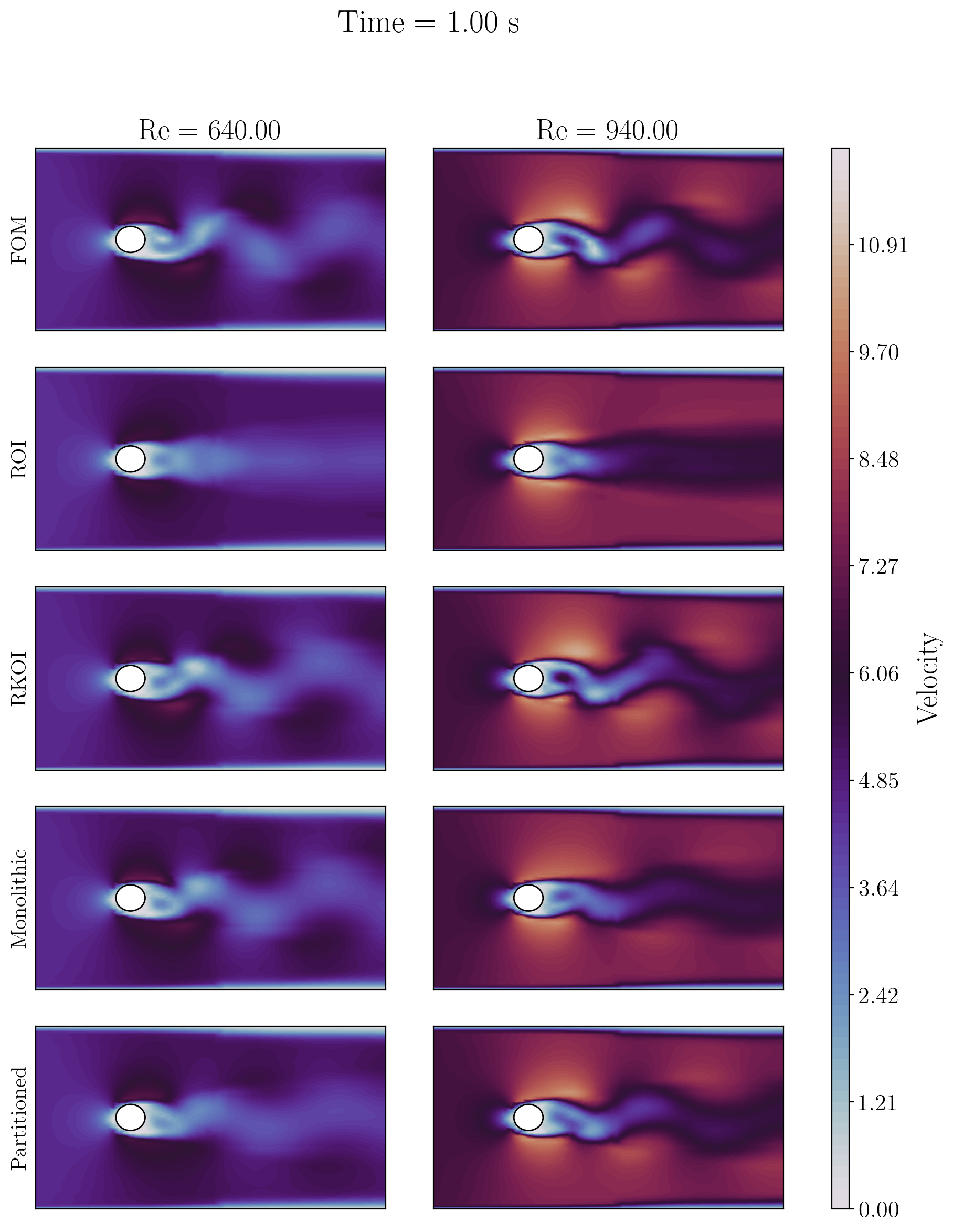}
    \caption{Contour plots of the velocity magnitude for the test parameters at the final time of the simulation, comparing the pDMD versions with the full-order solution (FOM).}
    \label{fig: flow-cfdbench-contour}
\end{figure}

The evolution of the error computed with Eq. \eqref{eqn: rec-time-error} over time is plotted in Figure \ref{fig: flow-cfdbench-temporal-errors}. The ROI approach shows its limitation when more complex physics comes into play: in fact, it is able to learn accurately the mean behaviour but it struggles to catch low-scale features (as seen in Figure \ref{fig: flow-cfdbench-contour}); in particular, the ROI works better compared to the other methods in predicting the initial transient. A possible remedy to handle high-frequency data would be to increase the rank (chosen to be 25 in this case): however, this would actually worsen the parametric interpolation \cite{faraji_data-driven_2024-2, NUTHOS24_pDMD}. The Monolithic (with rank 35) and the Partitioned (with rank 40) approaches show similar behaviour, being less able to predict accurately the establishment of the periodic flow, especially between 0.3 and 0.7 seconds. On the other hand, the RKOI approach (with rank 10) seems the most suited to handle complex flow physics due to its lower error compared to ROI and lower computational time compared to Monolithic and Partitioned, at least when the flow is fully developed from 0.2 s on. This latter aspect is of primary importance since the implementation of the latter two versions in the pyDMD package suffers from having a heavier online phase due to the fitting of the regressors within the prediction stage. The ROI and RKOI approaches, instead, have been designed to handle this part during the offline/training phase.

The test snapshots and the reconstruction with the different DMD approaches are plotted at the final instant in Figure \ref{fig: flow-cfdbench-contour}. As said above, the ROI approach is less able to find a proper surrogate model for the periodic flow and predicts mainly only the mean flow; conversely, the optimised version of the DMD, onto which the RKOI is based, is in general more suited and powerful in predicting oscillatory behaviours, hence better results are expected \cite{faraji_data-driven_2024-1}. The other algorithms are even more capable of doing so, even though the Monolithic and the Partitioned approaches require quite a high rank to do that, and the latter also requires some preprocessing of the data; on the contrary, the RKOI approach is more capable of building a low-rank approximation with less computational resources, keeping the online phase computationally efficient.

\subsection{The DYNASTY Experimental Facility}

The previous sections have been devoted to studying parametric DMD approaches for fluid dynamics data for the well-known thermal-hydraulic benchmark test case of the 'flow over cylinder' at different Reynolds numbers. In this last part of the discussion, these concepts will be applied to an experimental facility, DYNASTY \cite{cammi_dynasty_2016, benzoni_preliminary_2023}, deployed at Politecnico di Milano to study the dynamics of natural circulation for internally heated flows.

\subsubsection{Experimental Setup and RELAP5 Model}

\begin{figure}[tp]
    \centering
    \includegraphics[width=1\linewidth]{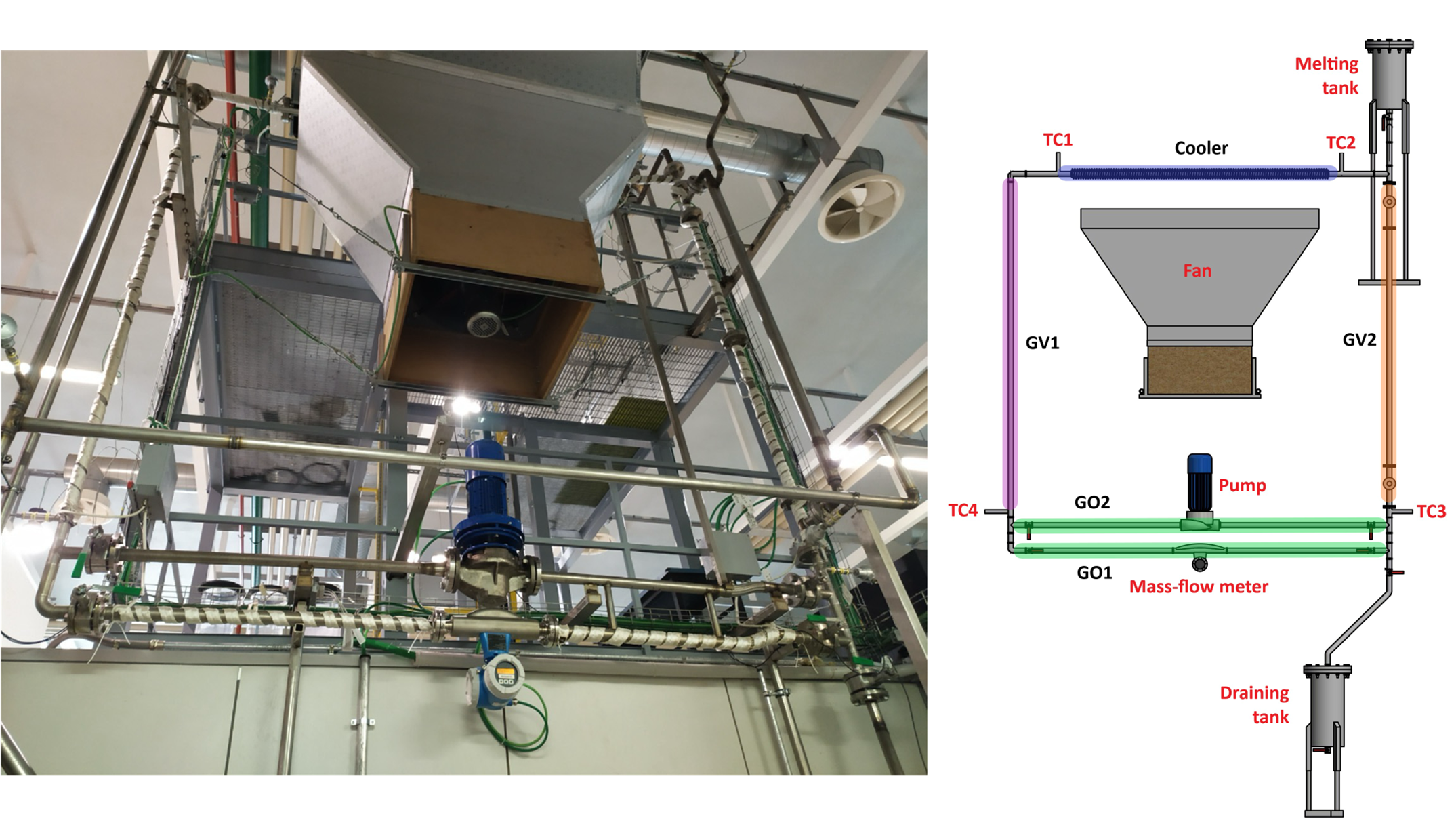}
    \caption{DYNASTY natural circulation loop \cite{cammi2019dynamics}: on the left, a picture of the real facility is provided, whereas on the right the scheme of the system is reported with the main components of the facility and the location of the heating strips.}
    \label{fig: dynasty-facility}
\end{figure}

The DYNASTY facility is a 3-by-3 square loop operating with a variety of fluids, mainly water and propylene glycol. The pipes are made of AISI316 stainless steel, with a diameter of 38 mm and a thickness of 2 mm. Out of the four legs of the facility, three of them are wrapped by electrical heating strips acting as independent heat sources \cite{cammi_dynasty_2016, cammi2019dynamics}\footnote{As the axial characteristic dimension is much larger than the radial one, internal heat generation can be approximated with distributed external heat generation}, as seen in Figure \ref{fig: dynasty-facility}. The upper horizontal tube is a finned tube acting as a heat sink, which can be cooled from below through an air fan. Various heating configurations can be experimentally set up, enabling either simultaneous activation of all strips or individual operation of one leg at a time. In the present work, a RELAP5/MOD3.3 \cite{fletcher1992relap5} model of DYNASTY was realised to simulate the experiment in the so-called VHHC-GV1 configuration (vertically-heated-horizontally-cooled, with GV1 as heat source) \cite{cammi_dynasty_2016, benzoni_preliminary_2023}.

The RELAP5 (R5) code was developed at the Idaho National Engineering Laboratory for the U.S Nuclear Regulatory Commission \cite{fletcher1992relap5, mangal2012capability} for light water reactor transient analysis. The code fulfils several functions, including licensing computations, evaluation of accident mitigation strategies, and evaluation of operator guidelines. However, R5 is a highly generic code that can be used to model a variety of thermal-hydraulic transients in both nuclear and non-nuclear systems involving mixtures of steam, water, non-condensable, and solute. Based on a non-homogeneous and non-equilibrium one-dimensional model for transient two-phase systems (which can also include non-condensable and soluble components) and a partially implicit numerical scheme, the goal of R5 is to include all relevant first-order effects for accurate prediction of transients whilst remaining sufficiently simple and cost-effective to allow sensitivity and parametric analyses. 

\begin{figure}[tp]
    \centering
    \includegraphics[width=0.75\linewidth]{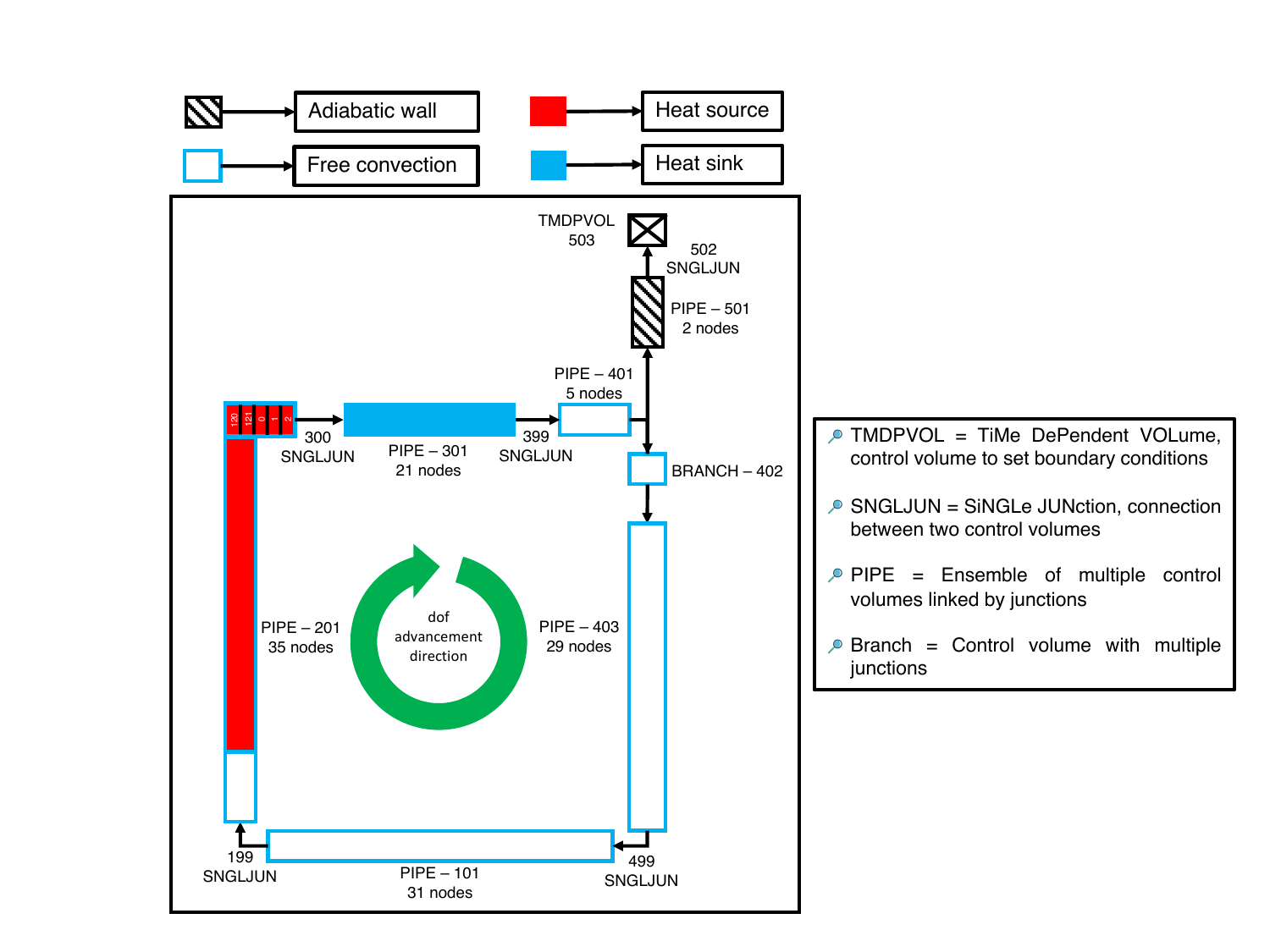}
    \caption{R5 nodalization of the DYNASTY experimental facility. The red zone corresponds to the heated section for the VHHC-GV1 case, whereas the blue zone corresponds to the finned cooler.}
    \label{fig: relap-nodes}
\end{figure}

The R5 nodalization of DYNASTY is depicted in Figure \ref{fig: relap-nodes}. The four DYNASTY legs are modelled as PIPE components, with the heater and cooler disposed in such a way as to reproduce the experimental configuration. Except for the cooler, where forced convection is imposed due to the presence of an airflow rate from the fan, all the remaining pipes exchange heat with the ambient, set to a temperature of 25 $^o$C, through free convection, for which the Churchill-Chu correlation was selected \cite{bergman2011fundamentals}. All the control volumes composing the loop are 100 mm long. The natural circulation is established from the left leg (PIPE-201), where, through an R5 heat structure, the fluid is heated, imposing the power level as in the related experiment \cite{cammi2019dynamics}. At the same time, the fan is switched on. The power level is enough to create the required buoyancy force that allows the fluid to overcome the gravitational and friction pressure drops. Regarding the latter, the distributed pressure drops were accounted for by considering a reasonable pipe roughness of 50 $\mu$m. Instead, the coefficients of the concentrated pressure drops for elbows and the T-junction at the top right corner were calculated from \cite{idelchik1987handbook}, while the pressure drop coefficient of the Coriolis flow meter was derived from its technical data sheet. The numerical simulation reproduced the startup of the natural circulation in the selected experiment. The transient lasted 2000 seconds, and the sampling timestep for storing the numerical data was 10 seconds \footnote{The first 400 seconds are not considered for DMD because the initial heating transient of the facility is neglected. In these first seconds, the temperature starts rising and a mass flow rate is established: this first part can be considered as a stand-alone transient with quite different dynamics and characteristic times compared to the rest of the transient, and the same algorithm can be applied to the snapshots pertaining to this subset \cite{NUTHOS24_pDMD}.}. The R5 model has already been validated against experimental data from different transients \cite{missaglia2025_dynasty_submitted}.

Regarding the experimental setup, four fluid thermocouples are placed in the middle of the channel roughly at the four corners of the facility (Figure \ref{fig: dynasty-facility}); experimental measurements have an uncertainty equal to $\pm$2 K. For more information on the experimental setup, refer to \cite{cammi_dynasty_2016, benzoni_preliminary_2023}.

\subsubsection{Application of parametric DMD}

The snapshots for this test case are the nodal fluid temperatures $T$ generated with the R5 code, considering the power $P$ (W) provided to each control volume as the parameter, with range $[22, 41]$ W. The dataset for this case has size $N_p = 20 \times \mathcal{N}_h = 122 \times N_t = 161$. The time interval considered to build the models is $[1400,3000]$ seconds\footnote{The time $t$ starts from 1400 seconds because of two reasons: from 0 to 1000 seconds, no heating is provided; power is turned on at $t$ = 1000 seconds, but this initial 400-second-long transient is neglected. In these first seconds, the temperature starts rising, and a mass flow rate is established: this first part can be considered as a stand-alone transient with quite different dynamics and characteristic times compared to the rest of the transient, and the same algorithm can be applied to the snapshots pertaining to this subset. Moreover, since the most interesting dynamics are observed from 5 minutes on, this initial variation is neglected.}.

\begin{figure}[htbp]
    \centering
    \includegraphics[width=0.85\linewidth]{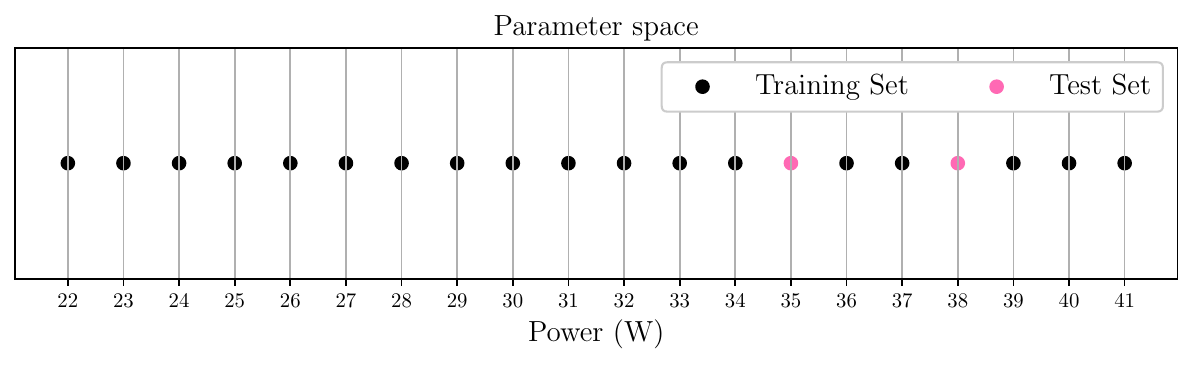}
    \caption{Train and test split of the power provided to each control volume for the R5 model of DYNASTY.}
    \label{fig: dynasty-traintest}
\end{figure}

The parameter set is divided into train and test sets randomly as in Figure \ref{fig: dynasty-traintest}, with the case $P=35$ W forced to be in the test set since experimental data for this case are public \cite{benzoni_preliminary_2023}. Compared to the previous cases, for this one, not only parameter interpolation but also the capabilities for time forecast are assessed: in fact, the time interval is divided into a train set, up until 2800 seconds, and a predict set spanning the whole time interval. The temperature and the experimental data have been rescaled to be between 0 and 1, allowing for their upload to the complementary GitHub repository. 

The SVD is always used to encode the spatial information and learn the surrogate dynamics in the latent space. In terms of ranks, following hyperparameter tuning, the ROI has been tested with rank $r=8$, the RKOI with rank $r=5$, and the Monolithic\footnote{As highlighted above, the Monolithic is more robust, providing lower errors with similar computational costs compared to the Partitioned version, available within pyDMD.} with rank $r=5$ and the Partitioned with rank $r=10$. 

 \begin{figure}[htbp]
    \centering
    \includegraphics[width=1\linewidth]{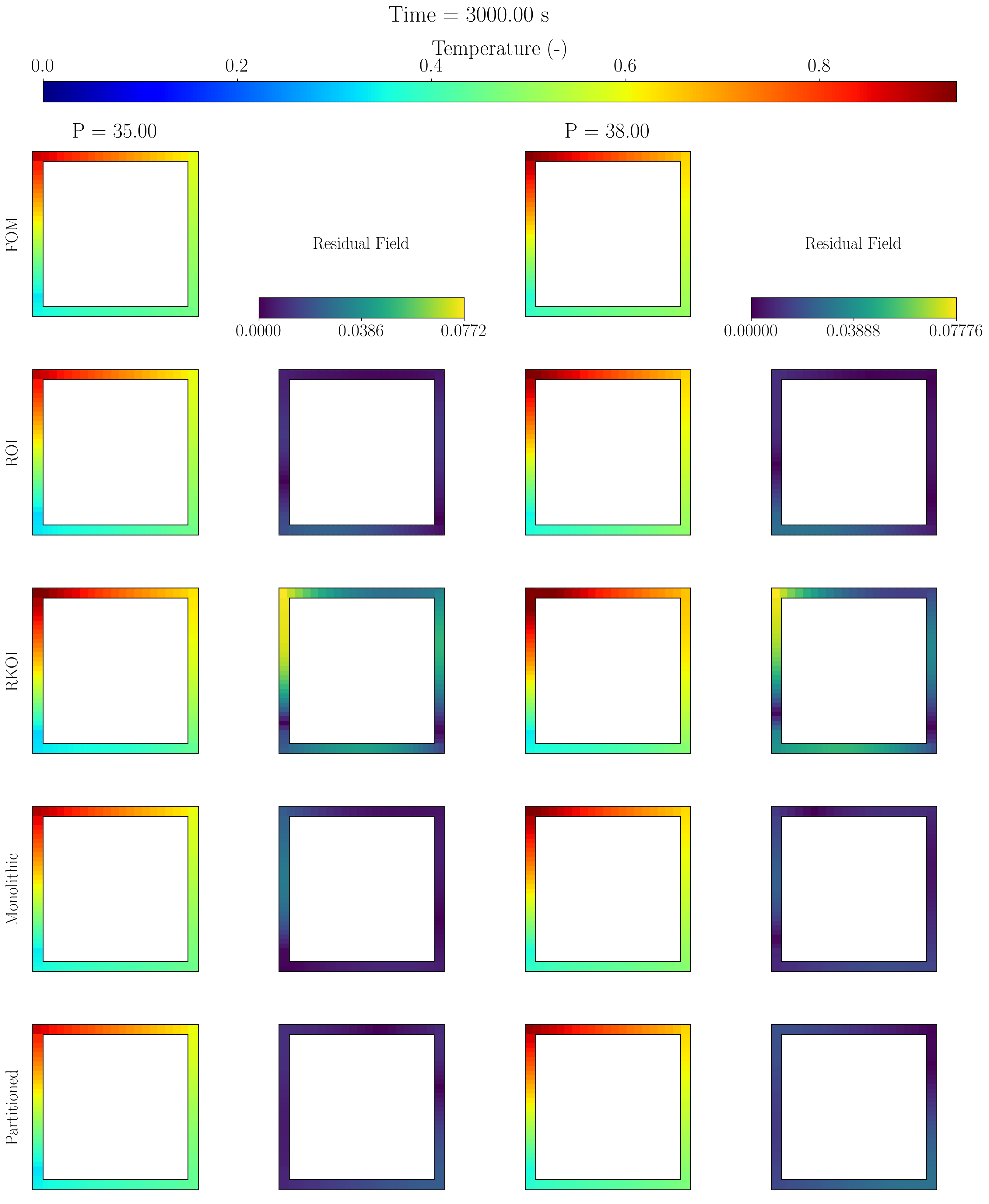}
    \caption{Contour plots of the temperature (scaled) for the test parameters forecasted at the final time of the simulation, comparing the pDMD versions with the full-order solution (FOM). In addition, the residual fields are shows in the even columns.}
    \label{fig: dynasty-contour}
\end{figure}

The temperature field for the different pDMD versions, along with the associated absolute residual field, is plotted in Figure \ref{fig: dynasty-contour}: each DMD version is able to predict a good reconstruction of the field, very similar to the FOM itself. By observing the time evolution of the temperature field (see video in the Supplementary Materials), the DMD reconstruction is very accurate in the train region (until 2800 seconds), whereas most of the errors are concentrated in the prediction region, even though the overall reconstruction is quite good: this explains the small differences that can be observed at the bottom left corner in Figure \ref{fig: dynasty-contour}. Overall, the four algorithms show comparable results, with the RKOI being the least performing one.

\begin{table}[htbp]
    \centering
    \begin{tabular}{c|c|c|c}
        \hline
        Parameter $\boldsymbol{\mu}^\star$ & Algorithm $*$ & Error $\epsilon^{\boldsymbol{\mu}^\star}_*$ & CPU Time (s) \\ 
        \hline
        \multirow{4}{*}{35.0} 
        & ROI & 0.0219 & 0.0031 \\ 
        & RKOI & \cellcolor{red!20}0.0383 & \cellcolor{green!20}0.0021 \\ 
        & Monolithic & \cellcolor{green!20}0.0112 & \cellcolor{red!20}3.2902 \\ 
        & Partitioned & 0.0223 & 0.1963 \\ 
        \hline
        \multirow{4}{*}{38.0} 
        & ROI & 0.0250 & \cellcolor{green!20}0.0206 \\ 
        & RKOI & \cellcolor{red!20}0.0372 & 0.0288 \\ 
        & Monolithic & \cellcolor{green!20}0.0208 & \cellcolor{red!20}3.4412 \\ 
        & Partitioned & 0.0311 & 0.2598 \\ 
        \hline
    \end{tabular}
    \caption{Relative error in the Frobenius norm and CPU times for different pDMD algorithms for two test parameters, using optimal ranks for the different versions for DYNASTY.}
    \label{tab: flow-dunasy-residuals-time-opt-rank}
\end{table}

Table \ref{tab: flow-dunasy-residuals-time-opt-rank} shows the relative error in the Frobenius norm with respect to the full order model, to retrieve an immediate comparison between the four different algorithms. The monolithic algorithm is indeed the most accurate, but it is intrinsically much heavier than all the other methods; on the other hand, ROI shows good performance with very low computation times required to obtain an approximation of the solution. 

\begin{figure}[htbp]
    \centering
    \includegraphics[width=1\linewidth]{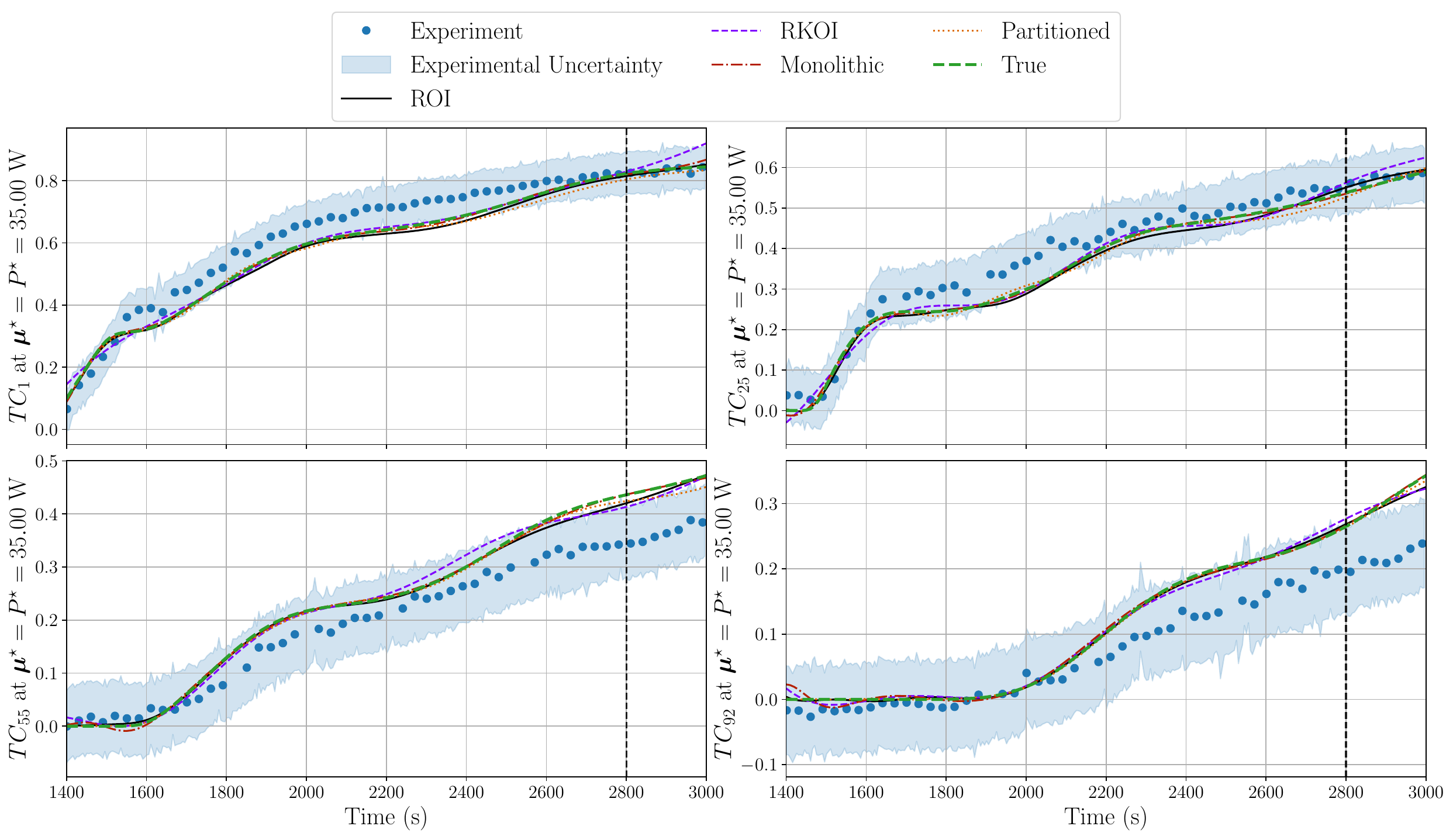}
    \caption{Validation of the R5 model and the pDMD approaches at the locations of the experimental thermocouples.}
    \label{fig: exp-dynasty-tc}
\end{figure}

In the end, the R5 solution and the pDMD reconstruction for the considered versions are compared with experimental data from \cite{benzoni_preliminary_2023}, as shown in Figure \ref{fig: exp-dynasty-tc}. It can be observed that the pDMD reconstruction, for all algorithm versions, is as accurate as the R5 solution, with some discrepancy in the prediction region in which there is still room for improvement. Nevertheless, both the FOM and the pDMD versions agree with the experimental data showing their accuracy, reaching an overall root-mean square error of 0.044 for ROI, 0.055 for RKOI and 0.049 for Monolithic in the prediction region.

\section{Conclusion} \label{sec:conclusion}

In this work, different versions of the parametric Dynamic Mode Decomposition algorithm for the prediction of parametric dynamical systems have been investigated and compared. The methods were applied to three different test cases: a laminar flow over cylinder, a benchmark dataset for the flow over cylinder, and the DYNASTY experimental facility. In all these scenarios, each pDMD version was tested. 

The DMD algorithm was originally developed to discover dynamical systems from temporal data; this work compares the different extensions of DMD currently available in the literature to handle parametric datasets. Three different versions have been tested: the Reduced Operator Interpolation, the Reduced Koopman Operator Interpolation and the Interpolation of the Latent Dynamics. Each one comes with advantages and shortcomings: the first two are very similar in their logic, and they are, in general, much cheaper from the computational point of view compared to the last one since they do not require to advance in time all the DMDs from the training data to interpolate the latent space; between the two, the RKOI is typically more performing for more chaotic dynamics since it relies on the optimised version of the DMD. The last one, in principle, is more accurate but comes with much higher computational costs as the interpolation step occurs in the online phase. 

The first result, common for all methods, is the requirement of a careful tuning of the rank $k$ in order to retrieve the optimal performance, and the optimal value of $r$ depends on the application under consideration. In terms of the comparison between algorithms, results show that the Reduced Operators Interpolation (ROI) approach is computationally efficient but it may suffer in predicting complex flow dynamics; conversely, the Reduced Koopman Operator Interpolation (RKOI) approach, by using optimised DMD instead of basic DMD, provides accurate and computationally efficient predictions; the Monolithic and Partitioned approaches from the pyDMD package also show good performances, but are more computationally expensive due to the fitting of regressors during the online phase. Although the RKOI method seems to be overall better, it must be stressed that this approach, compared to the others, leverages the excellent performances of optimised DMD. Therefore, generally speaking, the choice of the pDMD version depends on the specific requirements of the application, such as the complexity of the dynamics, the computational resources available, and the desired accuracy: in fact, when applied to a physical system such as the DYNASTY facility, the different pDMD versions show comparable performances in terms of parameter interpolation and time forecast. 

Overall, there does not seem to be one method unequivocally superior to the others, as each version requires a certain degree of tuning with which comparable performances can be obtained, at least in terms of accuracy. Additionally, this tuning strongly depends on the application and on the available dataset, and tuning of the hyperparameters for parametric test cases remains done on a trial-and-error basis. As a rule of thumb, when computational speed is sought, ROI and RKOI versions should be preferred; in terms of available documentation, the pyDMD package is a very valuable entry point for new users, especially if interested in DMD applications as a whole. However, when analysing a new dataset, it is worth trying all the versions and selecting the one most suited for that particular application and based on the objectives of the study.

In the future, the parametric DMD algorithms will be extended to data-assimilation problems to further improve their performance, including a surrogate model in a Kalman Filter framework to try to merge the background knowledge of the models with local evaluations given by the sensors, both during the training and the online phase. Additionally, the ROI algorithm will be upgraded, substituting basic DMD with optimised DMD, following the good performances obtained with the RKOI approach. Applications of the pDMD algorithms to multi-physics problems are also foreseen.

\section*{Code and Data Availability}
The code and data (compressed) that support the findings of this study are openly available at: \href{https://github.com/ERMETE-Lab/pDMD}{https://github.com/ERMETE-Lab/pDMD}.

\section*{List of Symbols}
\input{nomenclature}
\footnotesize{ \printnomenclature }

\newpage
\appendix

\section{Reduced Basis Methods: a brief introduction}\label{app: rb}
Within the Reduced Order Modelling framework, Reduced Basis (RB) methods play a crucial role, and they represent a fundamental brick in computational sciences. In the following, the rationale behind RB is briefly introduced: readers interested in the topic may refer to \cite{quarteroni_reduced_2015, rozza_model_2020, lassila_model_2014, brunton_data-driven_2022}.

Let $u(\vec{x}; t, \boldsymbol{\mu})$ be the solution of a mathematical model (e.g., a PDE), where $\vec{x}$ is the spatial coordinate, $t$ is the time, and $\boldsymbol{\mu}$ represents a vector of parameters, including input conditions (e.g., reactivity), physical parameters (e.g., viscosity) or boundary conditions (e.g., inlet velocity in a pipe). From this continuous function, a discrete vector state $\vec{u}(t_j, \boldsymbol{\mu}_i)\in\mathbb{R}^{\mathcal{N}_h}$ for time $t_i$ and for parameters $\boldsymbol{\mu}_j$ is saved, where $\mathcal{N}_h$ usually represents the dimension of the spatial mesh; this vector is typically called \textit{snapshot}, and a collection of them, for different times and parameters, represents the starting dataset needed to build any reduced model. The basic principle adopted in RB methods consists of the possibility of expressing the state (either in a continuous or in a discrete sense), adopting the following separation of variables:
\begin{equation}
    u(\vec{x}; t, \boldsymbol{\mu}) \simeq \sum_{i=1}^N \alpha_i(t, \boldsymbol{\mu}) \cdot \varphi_i(\vec{x})\qquad \text{ or }\qquad
    \vec{u}(t, \boldsymbol{\mu}) \simeq \sum_{i=1}^N \alpha_i(t, \boldsymbol{\mu}) \cdot \boldsymbol{\varphi}_i
    \label{eqn: chap02-rom-linear-exp}
\end{equation}
where $\varphi_i(\vec{x})/\boldsymbol{\varphi}_i$ are the basis functions representing the dominant spatial behaviour, whereas the coefficients $\alpha_i(t, \boldsymbol{\mu})$ encode the time and parametric behaviours. This decomposition allows for an approximation of the state once the reduced coefficients are known. It is important to mention that the number of basis functions $N$ is much lower than the spatial dimension $\mathcal{N}_h$, hence the dimensionality reduction. Furthermore, this compression operation can also be seen as a coordinate transformation from vectors of dimension $\mathcal{N}_h$ to a space of dimension $N$, where the basis functions represent the coordinate change. Then, the unknown coefficients $\alpha_i(t, \boldsymbol{\mu})$ can be computed in a computationally efficient way in this reduced space, thus retrieving a reduced representation of the time dynamics and parametric behaviour of the starting dataset. Different RB methods vary by how they compute the unknown coefficients. It is worth mentioning that there exists a theoretical measure of how well a starting full-order dataset can be represented by a reduced basis: common algebraic techniques such as the Singular Value Decomposition or the Proper Orthogonal Decomposition provide a practical measure of this quantity, thus allowing selecting an optimal value of $N$ \cite{Bachmayr2016}.

\section{Schemes of the parametric DMD models}\label{app: schemes-pdmd}

\subsection{Reduced Operators Interpolation}

\begin{figure}[htbp]
    \centering
    \includegraphics[width=1\linewidth]{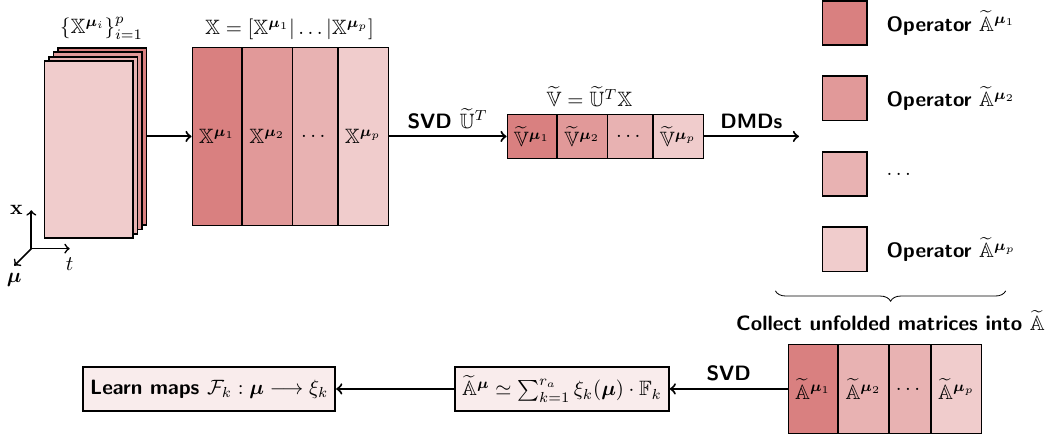}
    \caption{Scheme of the Reduced Operators Interpolation algorithm for parametric DMD (offline): the snapshots matrices $\mathbb{X}^{\boldsymbol{\mu}_i}$ are stacked together, and the spatial SVD is performed to obtain the parametric reduced dynamics $\widetilde{\mathbb{V}}^{\boldsymbol{\mu}_i}$; for each of them, a basic DMD is created to compute the DMD operators $\widetilde{\mathbb{A}}^{\boldsymbol{\mu}_i}$. These operators are unfolded and collected into a matrix $\widetilde{\mathbb{A}}$, then decomposed through a second SVD to encode the parametric dependence of the operators in a new set of reduced coefficients $\xi_k(\boldsymbol{\mu}$), onto which a regression model $\mathcal{F}_k$ can be built.}
    \label{fig: roi-offline}
\end{figure}

\begin{figure}[!htbp]
    \centering
    \includegraphics[width=1\linewidth]{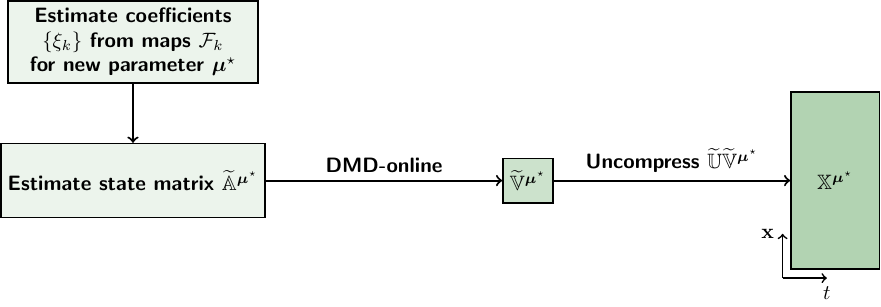}
    \caption{Scheme of the Reduced Operators Interpolation algorithm for parametric DMD (online): for a new parameter $\boldsymbol{\mu}^\star$ not included in the training set, the dynamical coefficients $\{\xi_k\}$ and the state matrix $\widetilde{\mathbb{A}}^{\boldsymbol{\mu}^\star}$ is estimated through the map $\mathcal{F}_k$ to solve the correspondent linear system; from it, the modal coefficients $\widetilde{\mathbb{V}}^{\boldsymbol{\mu}^\star}$ for the new parametric instance can be retrieved, then by decoding back onto the FOM space using the (first) SVD spatial modes an approximation of the high-dimensional field for the unseen parameter $\boldsymbol{\mu}^\star$ can be retrieved.}
    \label{fig: roi-online}
\end{figure}

\newpage
\subsection{Reduced Koopman Operators Interpolation}

\begin{figure}[!htbp]
    \centering
    \includegraphics[width=1\linewidth]{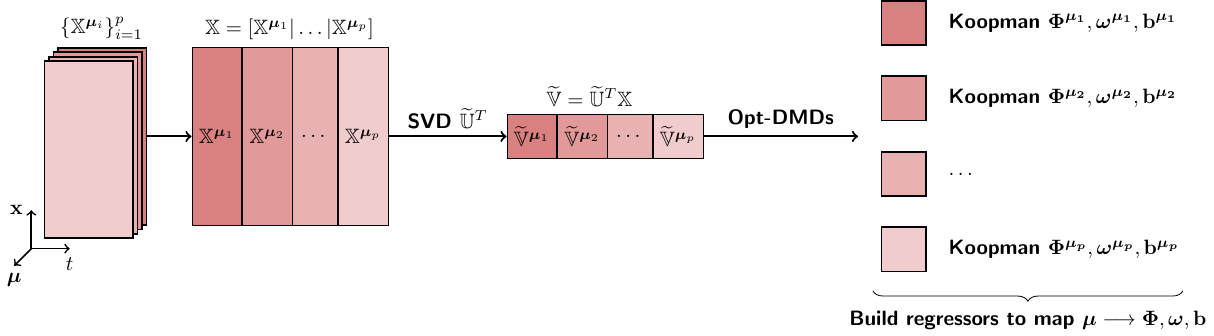}
    \caption{Scheme of the Reduced Koopman Operators Interpolation algorithm for parametric DMD (offline): the snapshots matrices $\mathbb{X}^{\boldsymbol{\mu}_i}$ are stacked together, and the spatial SVD is performed and obtained the parametric reduced dynamics $\widetilde{\mathbb{V}}^{\boldsymbol{\mu}_i}$; for each of them, an optimised DMD is computed and hence a collection of DMD modes $\boldsymbol{\phi}^{\boldsymbol{\mu}_i}$, amplitudes $\boldsymbol{b}^{\boldsymbol{\mu}_i}$ and frequencies $\boldsymbol{\omega}^{\boldsymbol{\mu}_i}$ is calculated. These are used to generate regressors mapping the parametric dependence.}
    \label{fig: rkoi-offline}
\end{figure}

\begin{figure}[!htbp]
    \centering
    \includegraphics[width=1\linewidth]{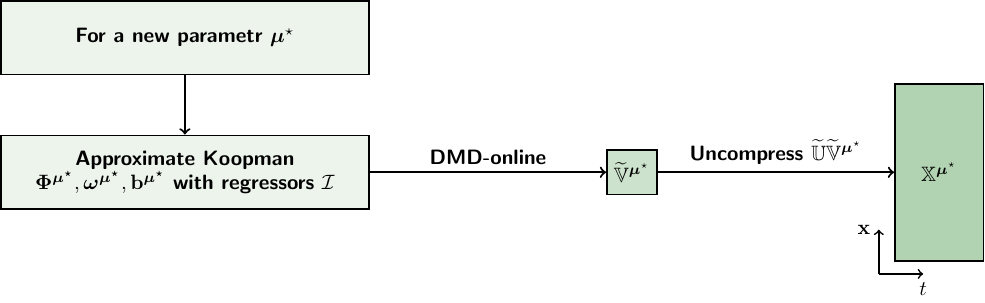}
    \caption{Scheme of the Reduced Koopman Operators Interpolation algorithm for parametric DMD (online): for a new parameter $\boldsymbol{\mu}^\star$, the dynamics is approximated by evaluating the regressors $\mathcal{I}$ and obtain the DMD modes $\boldsymbol{\phi}^{\boldsymbol{\mu}^\star}$, amplitudes $\boldsymbol{b}^{\boldsymbol{\mu}^\star}$ and frequencies $\boldsymbol{\omega}^{\boldsymbol{\mu}^\star}$, used to solve the correspondent linear system; the coefficients are uncompressed using the SVD spatial modes to retrieve the high-order approximation of the snapshot for the unseen parameter.}
    \label{fig: rkoi-online}
\end{figure}

\newpage
\subsection{Interpolation of Latent Dynamics}

\begin{figure}[!htbp]
    \centering
    \includegraphics[width=1\linewidth]{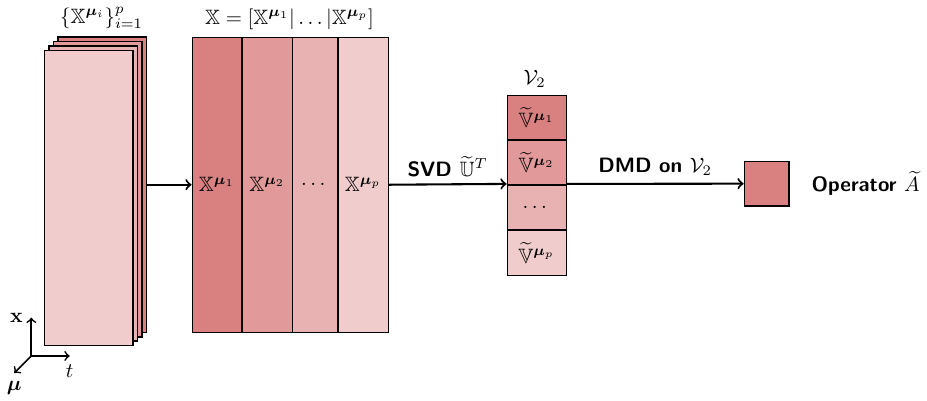}
    \caption{Scheme of the Monolithic approach from pyDMD (offline): the snapshots matrices $\mathbb{X}^{\boldsymbol{\mu}^i}$ are stacked together and the spatial SVD is performed and obtained the parametric reduced dynamics $\widetilde{\mathbb{V}}^{\boldsymbol{\mu}_i}$; the latent dynamics are stacked together in matrix $\mathcal{V}_2$ keeping the same columns (time instances) and a single DMD operator $\widetilde{A}$ is computed.}
    \label{fig: mono-offline}
\end{figure}

\begin{figure}[!htbp]
    \centering
    \includegraphics[width=1\linewidth]{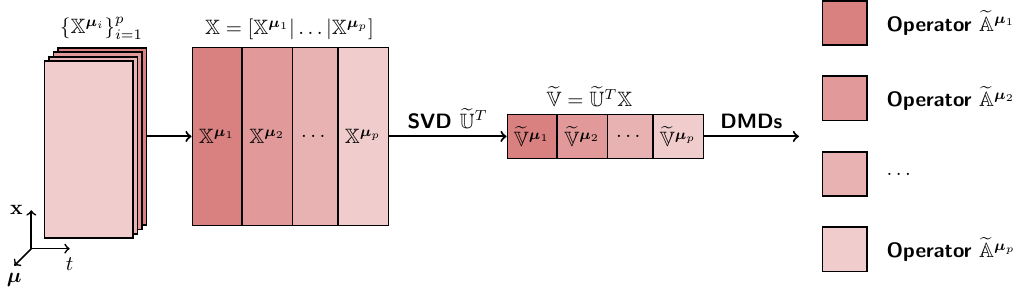}
    \caption{Scheme of the Partitioned approach from pyDMD (offline): the snapshots matrices $\mathbb{X}^{\boldsymbol{\mu}^i}$ are stacked together and the spatial SVD is performed to obtain the parametric reduced dynamics $\widetilde{\mathbb{V}}^{\boldsymbol{\mu}_i}$; for each of them, an optimised DMD model is created and the matrix $\widetilde{\mathbb{A}}^{\boldsymbol{\mu}_i}$ is computed}
    \label{fig: part-offline}
\end{figure}

\begin{figure}[!htbp]
    \centering
    \includegraphics[width=1\linewidth]{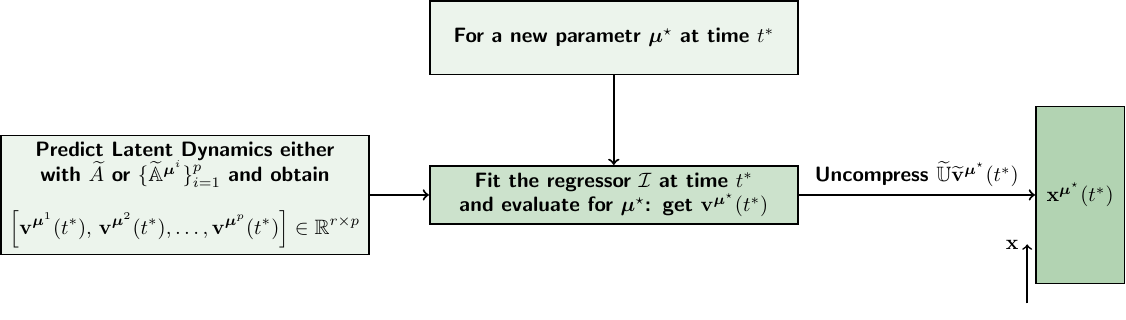}
    \caption{Scheme of the Monolithic and Partitioned approach from pyDMD (online): all the DMD models obtained during training are solved until the new time $t^\star$ to predict the latent dynamics $\boldsymbol{v}^{\boldsymbol{\mu}_i}(t^\star)$; for the new parameter $\boldsymbol{\mu}^\star$, a regressor $\mathcal{I}$ is fitted to predict the correspondent latent dynamics $\boldsymbol{v}^{\boldsymbol{\mu}^\star}(t^\star)$; then, the coefficients are uncompressed using the SVD spatial modes to retrieve the high-order approximation of the snapshot for the unseen parameter.}
    \label{fig: pydmd-online}
\end{figure}

\pagebreak
%Bibliography
\bibliographystyle{unsrt}  
\bibliography{bibliography}

\end{document}

%% file: nomenclature.tex
% Latin Symbols
\nomenclature[L]{$\vec{x}$}{High-dimensional state vector}
\nomenclature[L]{$\vec{u}$}{Velocity}
\nomenclature[L]{$\mathbb{X}$}{Snapshot matrix}
\nomenclature[L]{$\mathbb{A}$}{DMD operator}
\nomenclature[L]{$\mathcal{A}$}{Collection of DMD operators}
\nomenclature[L]{$\boldsymbol{a}$}{Flattened DMD operator}
\nomenclature[L]{$\mathbb{F}_k$}{$k-$th SVD mode of the DMD operators}
\nomenclature[L]{$\mathbb{W}$}{Matrix of DMD eigenvectors}
\nomenclature[L]{$\mathcal{F}_k$}{Map for the $k-$th reduced coefficient of the SVD on the DMD operators}
\nomenclature[L]{$\mathcal{I}$}{Regressor}
\nomenclature[L]{$\mathbb{U}$}{SVD spatial modes}
\nomenclature[L]{$\mathbb{V}$}{SVD reduced coefficients/reduced dynamics}
\nomenclature[L]{$\mathcal{V}_2$}{Snapshot matrix for the monolithic version of pDMD}
\nomenclature[L]{$\tilde{A}$}{Dynamical matrix for the monolithic version of pDMD}
\nomenclature[L]{$\vec{u}_j$}{$j$-th SVD spatial mode}
\nomenclature[L]{$\vec{v}_j$}{$j$-th reduced coefficient}
\nomenclature[L]{$\tilde{\mathbb{A}}$}{Reduced DMD operator}
\nomenclature[L]{$\mathcal{N}_h$}{Spatial degrees of freedom}
\nomenclature[L]{$\mathcal{D}$}{Parameter Space}
\nomenclature[L]{$\mathcal{X}$}{Parametric Snapshot matrix}
\nomenclature[L]{$p$}{Size of the Parameter Space, Pressure}
\nomenclature[L]{$N_t$}{Number of time instants}
\nomenclature[L]{$N_p$}{Number of parametric solutions}
\nomenclature[L]{$N_a$}{Rank of the SVD on the DMD operators}
\nomenclature[L]{$r$}{Rank of the SVD}
\nomenclature[L]{$\vec{b}$}{DMD amplitude}
\nomenclature[L]{$t$}{Time}
\nomenclature[L]{$T$}{Temperature}
\nomenclature[L]{$P$}{Power provided to control volume}
\nomenclature[L]{$Re$}{Reynolds number}

% Greek Symbols
\nomenclature[G]{$\boldsymbol{\mu}$}{Parameter}
\nomenclature[G]{$\boldsymbol{\Phi}$}{Matrix of DMD modes}
\nomenclature[G]{$\Sigma$}{Singular Values matrix}
\nomenclature[G]{$\epsilon_*^{\boldsymbol{\mu}}$}{Relative Error calculated with the Frobenius norm at parameter $\boldsymbol{\mu}$ for algorithm $*$}
\nomenclature[G]{$\varepsilon_*^{\boldsymbol{\mu}^\star}$}{Relative Error in time calculated with the Euclidian norm at parameter $\boldsymbol{\mu}$ for algorithm $*$}
\nomenclature[G]{$\nu$}{Kinematic Viscosity}
\nomenclature[G]{$\delta$}{Radius}
\nomenclature[G]{$\Lambda$}{Matrix of DMD eigenvalues}
\nomenclature[G]{$\xi_k$}{$k-$th reduced coefficient of the SVD on the DMD operators}
\nomenclature[G]{$\lambda_j$}{$j-$th DMD eigenvalue}
\nomenclature[G]{$\omega_j$}{Continuous $j-$th DMD eigenvalue}
\nomenclature[G]{$\boldsymbol{\phi}_j$}{$j-$th DMD eigenmode}

% Acronyms
\nomenclature[A]{DMD}{Dynamic Mode Decomposition}
\nomenclature[A]{pDMD}{Parametric Dynamic Mode Decomposition}
\nomenclature[A]{ROM}{Reduced Order Modelling}
\nomenclature[A]{PDE}{Partial Differential Equation}
\nomenclature[A]{SVD}{Singular Value Decomposition}
\nomenclature[A]{POD}{Proper Orthogonal Decomposition}
\nomenclature[A]{ROI}{Reduced Operator Interpolation}
\nomenclature[A]{RKOI}{Reduced Koopman Operator Interpolation}
\nomenclature[A]{FOM}{Full Order Model}
\nomenclature[A]{DYNASTY}{DYnamics of NAtural circulation for molten SalT internallY heated}
\nomenclature[A]{TC}{Thermo-Couple}
\nomenclature[A]{R5}{RELAP5}
\nomenclature[A]{VHHC}{Vertically Heated-Horizontally Cooled}